\documentclass[11pt]{amsart}
\usepackage{amssymb}
\usepackage{amsmath}
\usepackage{amsrefs}
\usepackage[colorlinks=true]{hyperref}
\usepackage{cleveref}
\usepackage{amscd}
\usepackage{amsthm}
\usepackage{mathtools}
\usepackage{tikz-cd}
\usepackage{mathrsfs}
\usepackage{pinlabel}
\usepackage[margin=2.93cm]{geometry}

\usepackage{xcolor}
\hypersetup{
    colorlinks,
    linkcolor={red!50!black},
    citecolor={green!50!black},
    urlcolor={blue!80!black}
}
\allowdisplaybreaks

\newtheorem{theorem}{Theorem}[section]
\newtheorem{lemma}[theorem]{Lemma}
\newtheorem{proposition}[theorem]{Proposition}
\newtheorem{corollary}[theorem]{Corollary}

\newtheorem{conjecture}[theorem]{Conjecture}
\newtheorem{question}[theorem]{Question}
\newtheorem{fact}[theorem]{Fact}

\newtheorem*{assumption-no-number}{Assumption}
\newtheorem*{corollary*}{Corollary}

\theoremstyle{definition}
\newtheorem{definition}[theorem]{Definition}
\newtheorem{example}[theorem]{Example}

\theoremstyle{remark}
\newtheorem{remark}[theorem]{Remark}

\numberwithin{equation}{section}

\newcommand{\R}{\mathbb{R}}

\newcommand{\Z}{\mathbb{Z}}

\renewcommand{\L}{\mathcal{L}}

\newcommand{\lt}{\left}
\newcommand{\rt}{\right}
\newcommand{\tms}{\times}

\newcommand{\rmk}{\begin{remark}}
\newcommand{\ermk}{\end{remark}}
\newcommand{\cor}{\begin{corollary}}
\newcommand{\ecor}{\end{corollary}}
\newcommand{\eq}{\begin{equation}}
\newcommand{\eeq}{\end{equation}}
\newcommand{\eqs}{\begin{equation*}}
\newcommand{\eeqs}{\end{equation*}}
\newcommand{\prop}{\begin{proposition}}
\newcommand{\eprop}{\end{proposition}}
\newcommand{\thm}{\begin{theorem}}
\newcommand{\ethm}{\end{theorem}}
\newcommand{\conj}{\begin{conjecture}}
\newcommand{\econj}{\end{conjecture}}
\newcommand{\lem}{\begin{lemma}}
\newcommand{\elem}{\end{lemma}}
\newcommand{\defi}{\begin{definition}}
\newcommand{\edefi}{\end{definition}}
\newcommand{\ex}{\begin{example}}
\newcommand{\eex}{\end{example}}
\newcommand{\alis}{\begin{align*}}
\newcommand{\ealis}{\end{align*}}
\newcommand{\pf}{\begin{proof}}
\newcommand{\epf}{\end{proof}}
\newcommand{\ali}{\begin{align}}
\newcommand{\eali}{\end{align}}
\newcommand{\qus}{\begin{question}}
\newcommand{\equs}{\end{question}}

\newcommand{\mc}{\mathcal}
\renewcommand{\bf}{\textbf}

\newcommand{\C}{\mathbb{C}}

\newcommand{\sub}{\subset}

\newcommand{\ov}{\overline}

\newcommand{\bb}{\mathbb}

\newcommand{\op}{\operatorname}

\renewcommand{\a}{\alpha}

\renewcommand{\d}{\partial}
\newcommand{\e}{\epsilon}

\newcommand{\g}{\gamma}

\newcommand{\s}{\sigma}
\renewcommand{\t}{\theta}

\renewcommand{\l}{\lambda}
\renewcommand{\o}{\omega}

\newcommand{\fk}{\frak}

\newcommand{\G}{\Gamma}
\renewcommand{\L}{\Lambda}

\renewcommand{\S}{\Sigma}

\renewcommand{\ov}{\overline}

\begin{document}
\title{Symplectic rigidity of fibers in cotangent bundles of open Riemann surfaces}
\author[Laurent C\^ot\'e]{Laurent C\^ot\'e}
\thanks {LC was supported by a Stanford University Benchmark Graduate Fellowship and by the National Science Foundation under Grant No.\ DMS-1926686.}
\address{Department of Mathematics, Harvard University, 1 Oxford St, Cambridge, MA 02138}
\email{lcote@math.harvard.edu}

\author[Georgios Dimitroglou Rizell]{Georgios Dimitroglou Rizell}
\thanks {GDR was supported by the grant KAW 2016.0198 from the Knut and Alice Wallenberg Foundation.}
\address {Department of Mathematics, Uppsala University, Box 480, SE-751 06, Uppsala, Sweden}
\email {georgios.dimitroglou@math.uu.se}

%\date{\today}  
\begin{abstract} We study symplectic rigidity phenomena for fibers in cotangent bundles of Riemann surfaces. Our main result can be seen as a generalization to open Riemann surfaces of arbitrary genus of work of Eliashberg and Polterovich on the Nearby Lagrangian Conjecture for $T^* \R^2$. As a corollary, we answer a strong version in dimension $2n=4$ of a question of Eliashberg about linking of Lagrangian disks in $T^* \R^n$, which was previously answered by Ekholm and Smith in dimensions $2n \geq 8$. \end{abstract}
\maketitle

\section{Introduction}

\subsection{Statement of results}

Let $\S$ be a connected, open Riemann surface of finite type and genus $g \geq 0$. We consider its cotangent bundle $(T^*\S, d\l_{\op{can}})$, where $\l_{\op{can}}$ is the canonical $1$-form. Let $F_x \sub (T^* \S, d \l_{can})$ be the cotangent fiber over some point $x \in \S$, which is a Lagrangian submanifold with respect to the symplectic form $d\l_{\op{can}}$. We prove the following theorem.  

\thm \label{theorem:main}
Let $L \sub T^* \S$ be a Lagrangian submanifold which is diffeomorphic to $\R^2$ and which agrees outside a compact set with the fiber $F_x$ for some $x \in \S$. Then $L$ is Hamiltonian isotopic to $F_x$ through a compactly supported Hamiltonian isotopy. 
\ethm

Recall that a symplectic isotopy is said to be Hamiltonian if it is generated by a (time-dependent) Hamiltonian vector field. Two Lagrangian submanifolds are said to be Hamiltonian isotopic if there is a Hamiltonian isotopy of the ambient symplectic manifold carrying one Lagrangian to the other. From the perspective of symplectic topology, it is natural to view Hamiltonian isotopic Lagrangian submanifolds as the same. 

The motivation for \Cref{theorem:main} is to exhibit symplectic rigidity phenomena for \emph{fibers} in cotangent bundles. Observe that cotangent bundles have two distinguished classes of Lagrangians, namely the zero section and the fibers. The celebrated Nearby Lagrangian Conjecture (often attributed to Arnold), which is one of the guiding problems of symplectic topology, asserts that any closed exact Lagrangian submanifold in a cotangent bundle is Hamiltonian isotopic to the zero section. This conjecture can be understood as a manifestation of symplectic rigidity for the zero section.  It is thus natural to search for manifestations of symplectic rigidity for cotangent fibers. 

It follows from work of Eliashberg--Polterovich \cite{eli-pol2} that any Lagrangian submanifold of $T^* \R^2$ which is diffeomorphic to $\R^2$ and agrees outside a compact set with a fiber is Hamiltonian isotopic to this fiber. In fact, due to the symmetry of $T^* \R^2$, this is easily seen to be equivalent to the Nearby Lagrangian Conjecture for $T^* \R^2$, which is the statement that Eliashberg and Polterovich originally proved. From this perspective, \Cref{theorem:main} can be viewed as a generalization of the work of Eliashberg--Polterovich. In particular, \Cref{theorem:main} recovers the Nearby Lagrangian Conjecture for $T^* \R^2$ as a special case.

One intriguing source of symplectic rigidity for cotangent fibers, which was originally promoted by Eliashberg (see \cite{ekholm-smith}), comes from studying linking phenomena. Observe that if $i: \R^n \to T^*\R^n$ is a Lagrangian embedding whose image $L:= i(\R^n)$ is disjoint from the fiber $F_0$ and agrees outside a compact set with $F_x$ for some $x \neq 0$, then we get a diagram
\eq
\begin{tikzcd}
\R^n \ar[r, "i"]  \ar[d, hook] & T^* \R^n - F_0 \ar[d, "\pi"] \\
S^n \ar[r] & \R^n - \{0\}.
\end{tikzcd}
\eeq 
Here $\pi$ is the projection, the left vertical arrow is the inclusion of $\R^n$ into its one-point compactification, and the bottom horizontal arrow represents the unique continuous map which makes the diagram commute. Eliashberg asked whether the map $S^n \to \R^n- \{0\}$ is nullhomotopic. This question was affirmatively answered by Ekholm and Smith in dimensions $2n \geq 8$; see \cite[Thm.\ 1.1]{ekholm-smith}.

Observe that Eliashberg's question is essentially asking whether $L$ can be homotopically linked with $F_0$. It is therefore a statement of \emph{homotopical} symplectic rigidity. The work of Ekholm and Smith can thus be seen as complementary to the recent developments establishing homotopical versions of the Nearby Lagrangian Conjecture (cf.\ \cite{abouzaid-kragh} and the references therein). 

One may ask whether similar homotopical rigidity results hold in more complicated examples. For instance, suppose that $M$ is an orientable (possibly closed) manifold of dimension $n$ and let $L \sub T^*M$ be a Lagrangian embedding of $\R^n$ which agrees outside a compact set with some fiber $F_p$. If $L$ is disjoint from another fiber $F_q \sub T^*M$, is the induced map $S^n \to M-q$ nullhomotopic? It seems conceivable that some of the tools which have been useful for studying homotopical versions of the Nearby Lagrangian Conjecture could also be applied to this type of question. 

By analogy with the ``full" Nearby Lagrangian Conjecture, it is natural to ask whether the above homotopical rigidity statements can be upgraded to the symplectic category. This leads to the following much stronger version of Eliashberg's original question.  

\begin{question} \label{question:full-fiber}
Suppose that $M$ is an orientable (possibly closed) manifold of dimension $n$ and let $L \sub T^*M$ be a Lagrangian embedding of $\R^n$ which agrees outside a compact set with some fiber $F_p$.  If $L$ is disjoint from another fiber $F_q \sub T^*M$, is $L$ Hamiltonian isotopic to $F_p$? If so, can the isotopy be confined to the complement of $F_q$? 
\end{question}

If \Cref{question:full-fiber} admits an affirmative answer, then a proof of this fact in all dimensions seems out of the reach of current technology. Indeed, there is a general lack of methods for constructing Hamiltonian isotopies between Lagrangians in symplectic manifolds of dimension six and higher. This is reflected in the fact that the Nearby Lagrangian Conjecture is not known for any cotangent bundle of dimension at least six.  In contrast, the situation is more favorable in dimension four where one can use certain pseudoholomorphic curve techniques which break down in higher dimensions.  In particular, the full Nearby Lagrangian Conjecture is known for some low genus cases, namely for $T^* \R^2, T^*(S^1 \tms \R), T^* S^2$ and $T^* \bb{T}^2$; see \cite{eli-pol2, dim1, hind, dgi}. 

The following corollary of \Cref{theorem:main} implies that \Cref{question:full-fiber} indeed admits an affirmative answer in dimension $2n=4$. 

\cor \label{corollary:main}
Let $\S$ be a (closed or open) Riemann surface of finite type and genus $g \geq 0$. Let $L$ be a Lagrangian submanifold which is diffeomorphic to $\R^2$ and which agrees outside a compact set with some cotangent fiber $F_x, \;x \in \S$. If $L \cap F_y= \emptyset$ for some $y \in \S, y \neq x$, then $L$ is isotopic to $F_x$ in the complement of $F_y$ though a compactly supported Hamiltonian isotopy. 
\ecor
Note that \Cref{corollary:main} can be immediately deduced from \Cref{theorem:main} by removing $F_y$ from $T^*\S$.

Taking a slightly different perspective on \Cref{corollary:main}, observe that it also implies the following: a Lagrangian $L$ which is diffeomorphic to $\R^2$ and which agrees with some fiber $F_x$ outside a compact set is Hamiltonian isotopic to $F_x$ if it can be displaced from a single other fiber. This hypothesis is of course necessary. For example, if we let $\tau(F_x)$ be the Dehn twist about the zero section of a cotangent fiber $F_x \sub T^* S^2$ for some $x \in S^2$, then $\tau(F_x)$ and $F_x$ are obviously not isotopic via a compactly supported Hamiltonian isotopy. It would be interesting to know whether the conclusion that $L$ is Hamiltonian isotopic to $F_x$ holds under the weaker assumption that $HF^{\bullet} (L, F_y)=0$.\footnote{According to a folklore strategy, if $HF(L, F_y)=0$ but $L \cap F_y \neq \emptyset$, one could attempt to disjoin $L$ from $F_y$ by using the holormorphic disks contributing to the Floer differential as ``Whitney disks". }

To the best of our knowledge, \Cref{corollary:main} is the first result which describes Lagrangian submanifolds up to Hamiltonian isotopy in cotangent bundles of closed Riemann surfaces of genus $g \geq 2$. In particular, the Nearby Lagrangian Conjecture is still open for such surfaces. Pseudoholomorphic curve techniques have been particularly powerful for studying Lagrangian submanifolds in symplectic $4$-manifolds. However, these techniques have proved to be difficult to apply in cases, such as cotangent bundles of Riemann surfaces of genus $g \geq 2$, where the symplectic manifold of interest does not come equipped with a natural foliation by pseudoholomorphic curves. Our arguments do not require such a foliation and are therefore more widely applicable. 
 
\subsection{Remarks on the proof of \Cref{theorem:main}}\label{subsection:remarks-proof-intro}

Let us briefly summarize our strategy for proving \Cref{theorem:main}.  The goal is to embed $L$ into a symplectic hypersurface $\tilde{Q} \sub T^*\S$ which diffeomorphic to $\R^3$. This hypersurface moreover admits a foliation with the following properties:
\begin{itemize}
\item $L$ is a leaf;
\item each leaf is (the image of) a proper Lagrangian embedding of $\R^2$;
\item each leaf coincides setwise with a cotangent fiber outside a compact subset of $T^*\S$ which is independent of the particular leaf. 
\end{itemize}

Having constructed $\tilde{Q}$, we deform $L$ through leaves of the foliation until it coincides with a cotangent fiber.  Let $(\chi_t)_{t \in [0,1]}$ denote this deformation, where $\chi_0=L$, each $\chi_t$ is Lagrangian and coinciding with a fibre outside of a compact subset, but where the deformation is not setwise fixed at infinity. To amend the latter defect, we construct an explicit family $(\s_t)_{t \in [0,1]}$ of compactly-supported diffeomorphisms of $\S$ so that $\s^*_t(\chi_t)$ is a family of Lagrangian submanifolds of $T^*\S$ which is fixed setwise outside a compact set. It is elementary to show that such a family is automatically induced by a Hamiltonian isotopy, which proves the theorem. 

To construct the hypersurface $\tilde{Q}$, we use the theory of punctured pseudoholomorphic curves. Roughly speaking, we start by constructing a hypersurface ``with holes" in \Cref{subsection:geometric-setup}. Then, in \Cref{subsection:filling-planes}, we build a hypersurface $Q$ by gluing moduli spaces of punctured pseudoholomorphic curves to fill the holes.  Finally, we obtain $\tilde{Q}$ by an explicit modification of $Q$ which is described in \Cref{subsection:completion-proof}. 

As mentioned previously, \Cref{theorem:main} generalizes a celebrated result of Eliashberg and Polterovich \cite{eli-pol2}.  Morally, our approach is inspired by their work -- in particular, the idea of embedding $L$ in a hypersurface is drawn from \cite{eli-pol2}. However, our implementation of this approach is substantially different, both in terms of the geometric content and of the holomorphic curve theory. As a result, our proof of \Cref{theorem:main} is logically independent of \cite{eli-pol2}, and does not reduce to \cite{eli-pol2} in the special case of $T^*\R^2$. 

To give an idea of how our methods differ from those in \cite{eli-pol2}, we remark that \cite{eli-pol2} relies on studying moduli spaces of pseudoholomorphic disks with boundary in a totally-real submanifold. In order to control the relevant moduli spaces, the authors introduce certain rather intricate geometric constructions which fundamentally depend on the fact that they are working in $\R^4=T^*\R^2$.  In contrast, our methods involve studying punctured pseudoholomorphic curves. Our geometric setup is therefore quite different, and appears a posteriori to be more flexible. 

We remark that the development of the theory of punctured pseudoholomorphic curves postdates \cite{eli-pol2}. This theory has proved extremely useful for studying Lagrangians in symplectic $4$-manifolds (see e.g.\ \cite{cieliebak-mohnke2, dgi, hind}), as well as in many other areas in symplectic and contact topology. 

\subsection*{Acknowledgements} We thank Yasha Eliashberg and Cliff Taubes for helpful conversations. We also wish to thank the anonymous referee for many helpful comments. Part of this work was carried out when the second author visited the Department of Mathematics at Stanford University in February 2019, and when the first author visited the Department of Mathematics at Uppsala University in December 2019.

\section{Preparations for the proof of \Cref{theorem:main}} \label{section:preparations}

In this section we collect preliminary material which is used in the arguments of \Cref{section:proof-main-thm}. More precisely, \Cref{subsection:punctured-curves} contains foundational material in the theory of punctured pseudoholomorphic curves. \Cref{subsection:a-c-structures} is devoted to constructing certain auxiliary almost-complex structures which are needed later. 

\subsection{Punctured pseudoholomorphic curves} \label{subsection:punctured-curves}

The theory of pseudoholomorphic curves in symplectic manifolds was initiated by Gromov \cite{gromov} and subsequently developed by many authors. While this theory was originally restricted to closed curves, or to curves with boundary in a Lagrangian or totally-real submanifold, much of it has subsequently been generalized to pseudoholomorphic maps from punctured Riemann surfaces into symplectic manifolds with cylindrical ends. This generalization plays an important role in many areas of symplectic topology (such as Symplectic Field Theory, and the study of low-dimensional contact and symplectic manifolds) and is used throughout this paper. 

For the reader's convenience, we collect in this section some foundational material in the theory of punctured pseudoholomorphic curves. Our presentation is entirely tailored to the needs of our paper and most of the definitions and results we state are special cases of more general statements. We refer the reader to \cite{wendlintersection} for a highly-readable introduction to the theory of punctured pseudoholomorphic curves. 

We begin with the following auxiliary definition.

\defi \label{definition:negative-end}
A manifold $W$ with a \emph{negative (or concave) cylindrical end modeled on $(Y, \l)$} consists in a datum $(W, Y, \l, e)$ where
\begin{itemize}
\item $W$ is a manifold of dimension $2n \geq 2$,
\item $Y$ is a manifold of dimension $2n-1$ equipped with a contact form $\l$,
\item $e$ is a proper embedding
\eqs 
e:  (-\infty, N] \tms Y \to W
\eeqs 
for some $N \in \R$.
\end{itemize}
We say that the negative end is non-degenerate (resp.\ Morse--Bott) if $\l$ is a non-degenerate contact form (resp.\ a Morse--Bott contact form, see \cite[Def.\ 1.7]{bourgeois-thesis}). By abuse of terminology, we will often refer to some manifold $W$ as having a negative cylindrical end modeled on $(Y, \l)$ without specifying the embedding $e$. 
\edefi

The notion of a manifold with a negative cylindrical end is only useful when we consider extra structure which is well-behaved with respect to the cylindrical end, as in the following definition.

\defi \label{definition:negative-end2}
Suppose that $W=(W, Y, \l , e)$ is a manifold with a negative cylindrical end modeled on $(Y, \l)$. 
	
Given an almost-complex structure $J$ on $W$, we say that $(W, J)$ is an almost-complex manifold with a negative cylindrical end if 
\eq
e^*(J)= \hat J_\l:= R_\l \otimes dt - \d_t \otimes \l + J_\l,
\eeq
where $J_\l: \xi \to \xi$ is an almost-complex structure which is compatible with $d\l|_{\xi}$, for $\xi:= \op{ker} \l$. 
	
Given a symplectic form $\o$ on $W$, we say that $(W, \o)$ is a symplectic manifold with a negative cylindrical end if 
\eq e^*\o= d(e^t \l).\eeq
\edefi

\rmk \label{remark:ends}
We do not require $Y$ to be compact or connected in the above definitions. It is however sometimes more convenient to talk about manifolds with negative \emph{ends}, where each end is assumed to be connected. One can also consider manifolds endowed with positive and negative cylindrical ends (see e.g. \cite[Sec.\ 3.2]{wendlintersection}), although we do not consider such structures in this paper. All of these notions are obvious adaptations of Definitions \ref{definition:negative-end} and \ref{definition:negative-end2}.
\ermk

The main examples which will be relevant in this paper are the following. 

\ex 
\label{ex:cotangentend}
Let $(L, g)$ be a Riemannian manifold and let $(T^*L, \l_{\op{can}})$ be the cotangent bundle of $L$, endowed with its canonical Liouville structure. The Liouville vector field is transverse to the sphere bundle $S^*_{\e, g}L:= \{\zeta \in T^*L \mid \|\zeta\|_g= \e\}$ for any $\e>0$. This implies that $(S^*_{\e, g}L, \l)$ is a contact manifold, where $\l$ denotes the restriction of $\l_{\op{can}}$.  It is well-known \cite[Sec.\ 1.5]{geiges} that the Reeb orbits are in bijective correspondence with the geodesics of $(L, g)$.

The Liouville flow furnishes a proper embedding 
\eq \label{equation:liouville-cotangent} \lt( (-\infty, 0] \tms S^*_{\e, g}L, d(e^t \l) \rt) \to (T^*L-L, \l_{\op{can}}),\eeq 
which makes $(T^*L-L, d\l_{\op{can}})$ into a symplectic manifold with a negative cylindrical end.
\eex

\ex
Let $(M, \o)$ be an arbitrary symplectic manifold and let $L \sub M$ be a Lagrangian submanifold. Fix a Weinstein embedding $\phi: \op{Op}(0_L) \to M$, where $\op{Op}(0_L) \sub T^*L$ is a neighborhood of the zero section. Given a Riemannian metric $g$ on $L$ and $\e>0$ small enough, we may precompose $\phi$ with the embedding \eqref{equation:liouville-cotangent}, thus endowing $(M-L, \o)$ the the structure of a symplectic manifold with a negative cylindrical end. Note that this structure of course depends on $\phi, g, \e$. 
\eex

\defi \label{definition:punctured-curve}
Let $(W,J)$ be an almost-complex manifold with (non-degenerate, Morse--Bott) cylindrical ends. A \emph{punctured pseudoholomorphic curve} is a map $u: \dot{\S} \to W$ satisfying the Cauchy-Riemann equations $$du \circ j = J \circ du.$$ Here $\dot{\S}= \S- \G$, where $\S$ is a compact Riemann surface (always assumed in this paper to be without boundary) and $\G \sub \S$ is a finite set of punctures. Two punctured pseudoholomorphic curves $(\S, j, \G,u), (\S', j', \G', u')$ are equivalent if there exists a holomorphic diffeomorphism $\phi: (\S, j, \G) \to (\S', j', \G')$ such that $u' \circ \phi= u$. 
\edefi

In general, punctured pseudoholomorphic curves can be rather badly behaved. However, in this paper, we will only consider punctured pseudoholomorphic curves which are \emph{asymptotically cylindrical}. This property means that the curve converges exponentially near each puncture to a (trivial cylinder over) a Reeb orbit. We refer to \cite[Sec.\ 1.1]{wendlintersection} for a precise definition, which will not be needed for our purposes. More generally, it is also useful to consider asymptotically cylindrical smooth maps. These are defined in the same way except that they are not required to be pseudoholomorphic. 

An asymptotically cylindrical punctured pseudoholomorphic curve with domain $\dot{\S}= \C$ will be referred to as a \emph{(pseudoholomorphic) plane}. 

Most foundational results in the theory of closed pseudoholomorphic curves in symplectic manifolds have been generalized to the punctured setting. We briefly mention some of these. 
\begin{itemize}
	\item The usual notion of energy when discussing punctured pseudoholomorphic curves is called the ``Hofer energy" (see \cite{hoferweinstein} and \cite[Sec.\ 1]{wendlsft}). The Hofer energy depends on the cylindrical structure of the ambient manifold. However, when the target manifold has no positive ends (as is always the case in this paper), the Hofer energy is controlled by the (ordinary) symplectic area \cite[Appendix]{cote}. 

	\item It can be shown that finite energy curves in symplectic manifolds with (non-degenerate, Morse--Bott) cylindrical ends are automatically asymptotically cylindrical (see \cite[Thm.\ 9.6]{wendlsft}). In fact, there are explicit formulas describing the asymptotic behavior of such curves (see \cite{hwz1, hwz4, siefring-asymptotic, siefring} as well as \cite[Sec.\ 3.2]{wendlintersection}) which play an important role in the theory. 
		
	\item The generalization of Gromov's compactness to the punctured setting is the so-called \emph{SFT compactness theorem} \cite{sftcompactness, cmsftcompactness}. If $(W,\o, J)$ is simultaneously a symplectic and almost-complex manifold with a negative (non-degenerate, Morse--Bott) end, and if $\o$ tames $J$ (see \cite[Def.\ 1.7]{wendlintersection}), then this theorem implies that a sequence of punctured curves with uniformly bounded symplectic area converges to a so-called \emph{pseudoholomorphic building} (a building is, roughly speaking, a finite collection of curves which satisfy various compatibility conditions \cite[Sec.\ 7-9]{sftcompactness}). 
	
	\item The functional-analytic setup for constructing moduli spaces of punctured pseudoholomorphic curves is described in \cite[Sec.\ 3.2]{wendlautomatic}. As usual, there is a notion of a curve being \emph{regular} or \emph{transversally cut out}, which essentially means that the nonlinear Cauchy-Riemann operator $\ov{\d}_J$, viewed as a section of an appropriate infinite dimensional vector bundle, intersects the zero section transversally at this curve. Since the details of the setup will not be relevant in this paper, we do not discuss them further here.  

\end{itemize}

In this paper, we will be considering pseudoholomorphic planes in an almost-complex manifold with a Morse--Bott negative cylindrical end.  When considering moduli spaces of such planes, one needs to specify whether the planes converge to a fixed Reeb orbit at the puncture, or whether we allow them to converge to any Reeb orbit in a Morse--Bott family. We speak respectively of a \emph{constrained} or \emph{unconstrained} puncture (cf.\ \cite[Def.\ 1.1]{wendlautomatic}).   

Unless otherwise specified, we always assume in this paper that the punctures are unconstrained. In fact, the only place in this paper where we need to consider curves with a constrained puncture occurs in Lemmas \ref{lemma:index-computation} and \ref{lemma:restricted-moduli}. 

Our next task is to state the index formula for pseudoholomorphic planes in almost-complex $4$-manifolds with a negative cylindrical end (this is of course a special case of a more general index formula).

To set the stage, let $(Y, \l)$ be a contact manifold with $\xi:= \op{ker} \l$. Let $J_\l: \xi \to \xi$ be an almost-complex structure which is compatible with the symplectic form $d\l$. Given a Reeb orbit $\g: S^1 \to Y$ parametrized so that $\l(\dot{\g})=1$, there is an associated \emph{asymptotic operator} 
\eq \bf{A}_\g:= \G(\g^*\xi) \to \G(\g^*\xi); \; \; \eta \mapsto -J_\l(\nabla_t \eta - \nabla_\eta R_\l).\eeq
Here $R_\l$ denotes the Reeb vector field associated to $\l$; $\nabla$ denotes the covariant derivative associated to some symmetric connection and $\nabla_t$ denotes covariant differentiation along $\g$. One can verify that the asymptotic operator does not depend on the choice of symmetric connection. However, it does depends on $J_\l$.  See e.g.\ \cite[Sec.\ 3.1]{wendlintersection} for more on asymptotic operators. 

Let us now fix a complex trivialization $\Phi$ of $\g^*\xi= (\g^*\xi, J_\l)$. The \emph{Conley--Zehnder index} 
\eq \mu_{CZ}^\Phi(A_\g)= \mu_{CZ}^\Phi(\g) \in \Z \eeq 
is an integer valued invariant of $\g$ and $\Phi$ (more precisely, it depends on the contact manifold $(Y, \l)$, the Reeb orbit $\g$ and the trivialization $\Phi$ up to homotopy through complex trivializations).  On a contact $3$-manifold, the Conley--Zehnder index essentially measures the amount of rotation of the Reeb flow along a given Reeb orbit. We refer to \cite[Sec.\ 3.4]{wendlsft} for more on the Conley--Zehnder index. 

Finally, we need to introduce the (relative) first Chern number.
\defi\label{definition:relative-chern}
Let $(W, J)$ be an almost-complex manifold of real dimension $2n$ with a (non-degenerate, Morse--Bott) negative cylindrical end modeled on $(Y, \l)$. Let $\Phi$ be a complex trivialization of $\xi= \op{ker} \l$ with respect to $J_\l:= J|_{\xi}$. By adjoining the Reeb vector field, we obtain a natural identification between complex trivializations $\Phi$ of the contact planes $(\xi, J_\l)$ and complex trivializations of $(TW, J)$ defined at the negative end.  Let $u : \S \to (W,J)$ be an asymptotically cylindrical punctured pseudoholomorphic curve. The \emph{(relative) Chern number} of $(u^*TM, J)$ is a signed count of zeros of a generic section of $\Lambda_\C^n(u^*TW)$ which is required to be non-vanishing and constant at infinity with respect to $\Phi$. We denote this quantity by 
\eq c_1^{\Phi}(u) \in \Z. \eeq
By a slight abuse of terminology, we sometimes refer to $c_1^{\Phi}(u)$ as the Chern number of $u$.
\edefi

We now come to the promised index formula. 

\defi \label{definition:index}
Let $(W,J)$ be an almost-complex $4$-manifold with a (non-degenerate, Morse--Bott) negative end modeled on $(Y, \l)$. Let $u: (\C, j) \to (W, J)$ be a pseudoholomorphic plane asymptotic to a Reeb orbit $\g$. 

The \emph{index} of $u$ is defined as follows (see \cite[(1.1)]{wendlautomatic}): 
\eq \label{equation:index} \op{ind}(u)= -1 + 2c_1^{\Phi}(u^*TW) - \mu_{\op{CZ}}^{\Phi}(\bf{A}_{\g}+\bf{c}), \eeq
 where $\bf{c}= \delta \cdot \op{Id}$ if the puncture is unconstrained and $\bf{c}= - \delta \cdot \op{Id}$ if the puncture is constrained, for $\delta>0$ small enough. 
Here $\Phi$ denotes a complex trivialization of $\g^*\xi$ with respect to $J_\l=J|_{\xi}$.  The index is independent of the choice of trivialization and of $\delta>0$ provided that $\delta$ is small enough. 
\edefi

\rmk \label{remark:index-difference}
Suppose that $u$ is a pseudoholomorphic plane whose asymptotic orbit $\g$ is contained in a Morse--Bott manifold $X_{\op{MB}}$ (see \cite[Sec.\ 1.1]{wendlautomatic}). Then the difference between the index in the unconstrained case and the index in the constrained case is precisely (see \cite[(3.3)]{wendlautomatic})
\eq \label{equation:index-difference} \mu_{\op{CZ}}^{\Phi}(\bf{A}_{\g}-\delta \cdot \op{id})-\mu_{\op{CZ}}^{\Phi}(\bf{A}_{\g}+\delta \cdot \op{id})= \op{dim}(X_{\op{MB}})-1. \eeq This matches our intuition that the unconstrained moduli space has additional degrees of freedom corresponding precisely to the dimension of the Morse--Bott family. (Note that $\op{dim}(X_{\op{MB}})=1$ if and only if $X_{\op{MB}}$ is an isolated Reeb orbit). 
\ermk

\defi \label{definition:normal-chern}
Let $(W,J)$ be an almost-complex manifold with a (non-degenerate, Morse--Bott) negative cylindrical end. Let $u: \dot{\S} \to (W, J)$ be an asymptotically cylindrical punctured pseudoholomorphic curve. The \emph{normal Chern number} $c_N(u) \in \frac{1}{2}\Z$ can be interpreted as the first Chern number of the normal bundle and is defined as follows (see \cite[(1.2)]{wendlautomatic}):
\eq \label{equation:normal-chern} 2c_N(u) = \op{ind}(u) - 2 +2g +\#\G_0. \eeq 

Here $\#\G_0$ is a count of \emph{even} punctures (see \cite[Sec.\ 3.2]{wendlautomatic}) which is either $0$ or $1$ if $u$ is a plane, and depends on whether the puncture is constrained or unconstrained.  
\edefi

\rmk \label{remark:homotopies}
The relative Chern number, index and the normal Chern number are in fact defined for asymptotically cylindrical \emph{smooth} map and are invariant under homotopies through such maps (cf.\ \cite[Sec.\ 3.4]{wendlintersection}). This can be verified by inspecting Definitions \ref{definition:relative-chern}, \ref{definition:index} and \ref{definition:normal-chern}, all of which involve purely topological quantities.
\ermk

The intersection theory for punctured pseudoholomorphic curves in almost-complex manifolds of dimension $4$ with cylindrical ends was developed by Siefring \cite{siefring} in the non-degenerate case, and by Siefring--Wendl \cite{siefringwendl} in the Morse--Bott case.\footnote{At the time of writing, this work is still in preparation. However, the fact that one can extend Siefring's intersection theory to the Morse--Bott setting is widely accepted by experts, and has already been used in many applications.} We refer to \cite{wendlintersection} and \cite[Appendix A.3]{wendlduke} for an overview of this theory in both the non-degenerate and Morse--Bott settings.

\begin{fact}[Siefring intersection number; see A.3 in \cite{wendlduke} and \cite{wendlintersection}] \label{fact:siefring-number} 
Suppose that $(W,J)$ is an almost-complex $4$-manifold with non-degenerate or Morse--Bott cylindrical ends. Let $u: \dot{\S} \to W, v: \dot{\S}' \to W$ be asymptotically cylindrical smooth maps. The \emph{Siefring intersection number} $$u*v \in \Z$$ satisfies the following properties:
\begin{enumerate}
\item \label{item:homotopy-invariance} $u*v$ is invariant under homotopies of $u$ and $v$ through asymptotically cylindrical smooth maps;
\item \label{item:disjointasymptotics} if $u, v$ are disjoint and do not share any asymptotic orbits, then $u*v=0$.
\end{enumerate}
\end{fact}

The Siefring intersection theory in almost-complex $4$-manifolds is particularly powerful when coupled with the following two results, which we refer to as the \emph{adjunction formula} and \emph{automatic transversality}. The adjunction formula was originally proved by Siefring \cite[Sec.\ 4.2]{siefring} (in the non-degenerate case) and the automatic transversality result is due to Wendl \cite{wendlautomatic}.

\begin{fact}[Adjunction formula; see A.3 in \cite{wendlduke}] \label{fact:adjunction} 
Suppose that $(W,J)$ is an almost-complex $4$-manifold with non-degenerate or Morse--Bott cylindrical ends. Let $u: (\dot{\S}, j) \to (W, J)$ be a somewhere injective, asymptotically cylindrical punctured pseudoholomorphic curve, where $\dot{\S}= \S- \G$. Then 
\eq \label{equation:adjunction} u*u= \op{sing}(u)+ 2c_N(u)+\sum_{z \in \G} \op{cov}_{\infty}(z).\eeq 
Each of the terms on the right-hand side of \eqref{equation:adjunction} is non-negative and invariant under homotopies of $u$ through asymptotically cylindrical smooth maps. The terms $\op{cov}_{\infty}(-)$ vanish in case the asymptotic orbits of $u$ are simply-covered. Moreover, $u$ is embedded if $\op{sing}(u)=0$. 
\end{fact}

\begin{fact}[Automatic transversality; see Thm.\ 1 in \cite{wendlautomatic}] \label{fact:automatic-transversality} 
Let $(W, J)$ be an almost-complex $4$-manifold with non-degenerate or Morse--Bott cylindrical ends. Suppose that $u: \S \to W$ is a non-constant, asymptotically cylindrical punctured pseudoholomorphic curve \emph{which is immersed}. If 
\eq\label{equation:aut-transversality}
\op{ind}(u) > c_N(u),
\eeq 
then $u$ is regular. 
\end{fact}

\Cref{fact:automatic-transversality} is a special case of a more general automatic transversality criterion due to Wendl (Thm.\ 1 in \cite{wendlautomatic}) which holds without the assumption that $u$ is immersed. (The more general criterion involves an additional term which vanishes if and only if $u$ is immersed.)

\subsection{Construction of some auxiliary almost-complex structures} \label{subsection:a-c-structures} In this section, we explicitly construct certain almost-complex structures which will be needed later on. The constructions are not very illuminating, so the reader may wish to skip directly to \Cref{subsection:geometric-setup} and return to this section when the need arises.

We write $T^*(S^1 \tms \R)= \R/(2\pi\Z) \tms \R^3$ with coordinates $(\t, t, r,s )$ and symplectic form $\o_{can}:= dt \wedge d\t + ds \wedge dr$. The zero section $0_{S^1 \tms \R}$ is given by $\{ t=s=0\}$. 

We consider constants $C_1>100$ and $\tilde{\e} < 1/100$ which will be fixed in \Cref{subsection:geometric-setup}.  

Let 
\eq\label{equation:sphere-bundle} S_{\tilde{\e}}:= \{ (\t, t, r,s) \mid \|(t,s) \| = \tilde{\e} \} \sub T^*(S^1 \tms \R),\eeq 
where the norm is induced by the standard flat metric on $S^1 \tms \R$. Letting $V= t \d_t + s \d_s$ denote the radial Liouville vector field, $S_{\tilde{\e}}$ is a contact manifold with respect to $\a:= i_V \o_{can}$. Letting $t= \tilde{\e} \cos \phi$ and $s = \tilde{\e} \sin \phi$ for $\phi \in \R/(2\pi \Z)$, we have natural coordinates $(\t, r, \phi )$ for $S_{\tilde{\e}}$ and we compute that $\a= \tilde{\e}( \cos \phi d\t + \sin \phi dr)$. 

For $u \in (-1/2, 1/2)$, let $f(u)= \sqrt{C_1^2+2u}$ and consider the embedding 
\begin{align}\label{align:def-of-F}
F: \{ t \in (-1/2, 1/2) \} \sub  \R/(2\pi\Z) \tms \R^3 &\to \left(\R^4_{(x_1,y_1, x_2, y_2)},dx_1 \wedge dy_1+dx_2 \wedge dy_2\right) \\
(\t, t, r, s) &\mapsto ( f(t) \cos (\t), f(t) \sin (\t), r, -s)
\end{align}
It is straightforward to check that $F$ is in fact a symplectic embedding. Letting $j$ denote the standard complex structure on $\R^4$, i.e.~$j(\partial_{x_i})=\partial_{y_i}$, one computes that 
\eq \label{equation:j-pullback}
F^*(j) = dF^{-1} \circ j \circ dF = 
 \begin{pmatrix}
 0 & 1/(C_1^2+2t) & 0 &0 \\
-(C_1^2 + 2t) & 0 &0 & 0 \\
 0 &0 & 0&1 \\
  0 & 0 & -1 & 0
 \end{pmatrix}. \eeq

Let $\mc{Q}:= \{(\t, t, r,s) \mid 0< \|(t,s)\| < 1\} \sub (\R/2 \pi\Z) \tms \R^3$. For $i=1,2$, fix smooth functions $\rho_i: \mc{Q} \to \R$ which satisfy the following properties:
\begin{itemize}
\item $\rho_i>0$, 
\item $\rho_i = \tilde{\e} \| (t, s) \| $ if $ \| (t, s) \| \leq \tilde{\e} $, 
\item $\rho_1= C_1^2 + 2 t $ if $\| (t, s) \| \geq 2 \tilde{\e}$,
\item $\rho_2= 1 $ if $\| (t, s) \| \geq 2 \tilde{\e}$.
\end{itemize}

It's clear that functions satisfying the above properties exist.  

We let $J_0$ be the unique almost-complex structure on $\mc{Q} \sub  (\R/2 \pi\Z) \tms \R^3- \{ t=s=0\}=T^*(S^1 \tms \R) - 0_{S^1 \tms \R}$ which satisfies $J_0(\d_\t)= -\rho_1 \d_t$ and $J_0(\d_r)= -\rho_2 \d_s$. 

\lem \label{lemma:cylindrical}
The almost-complex structure $J_0$ is cylindrical with respect to the canonical symplectic embedding  
\eq\label{equation:cyl-end} \iota: ((-\infty, 0] \tms S_{\tilde{\e}}, d(e^{\tau} \a) ) \hookrightarrow (T^*(S^1 \tms \R) - 0_{S^1 \tms \R}, \o_{can} ) \eeq 
induced by the Liouville flow, where $\tau$ is the variable corresponding to  $(-\infty, 0]$. 
\elem

\pf
By explicit computation, we find that $\iota^*(J_0)( \d_{\tau}) = \frac{1}{\tilde{\e}} ( \cos \phi \d_{\t} + \sin \phi \d_r)= R_{\a}$ and that $\iota^*(J_0)( \frac{1}{\tilde{\e}} (\sin \phi \d_{\t} - \cos \phi \d_r)) = \d_{\phi}$. Since $\op{ker} \a = \op{span} \{\sin \phi \d_{\t} - \cos \phi \d_r , \d_{\phi} \}$, this proves the claim. 
\epf

In the sequel, we let $J_{cyl}$ be the unique almost-complex structure on $\R_\tau \tms S_{\tilde{\e}}= T^*(S^1 \tms \R)- 0_{S^1 \tms \R}$ which satisfies $J_{cyl}(\d_{\tau})= R_{\a}$ and $J_{cyl}( \frac{1}{\tilde{\e}} ( \sin \phi \d_{\t} + \cos \phi \d_r)) = \d_{\phi}$. In particular, 
\eq\label{equation:jcyl-j0} J_{cyl}= \iota^*(J_0), \text{ for } \tau \leq 0. \eeq

We also record the following fact, which is an easy computation using the definition of the Conley--Zehnder index in \cite[Sec.\ 1.1]{wendlautomatic}. We do not give a proof since an essentially identical computation can be found in \cite[Appendix A]{hindandlisi}. (The relevant computation is stated on p.\ 1118 of \textit{loc}.\ \textit{cit}.\ and corresponds to the case of a negative puncture allowed to move in a Morse--Bott family of periodic Reeb orbits that are lifts of oriented closed geodesics on a flat two-torus; recall that the contact form considered here is induced by a flat metric on the cylinder.)

\lem \label{lemma:cz-index} 
Consider the complex trivialization $\Theta$ of $\op{ker} \a$ given by $\Theta=\{ \frac{1}{\tilde{\e}} ( \sin \phi \d_{\t} - \cos \phi \d_r)=-J_{cyl}\d_{\phi}\}$. Consider the Morse--Bott family of Reeb orbits $X_1= \{\phi=0\} \sub S_{\tilde{\e}}$ and let $\g \sub X_1$ be a closed Reeb orbit. 
For $\delta>0$ small enough, we have  
\eq \label{equation:cz-index} \mu_{\op{CZ}}^{\Theta}(\bf{A}_{\g}+\delta \cdot \op{id})=  0. \eeq
\qed 
\elem

We now define an almost-complex structure $\tilde{J}$ on $\R^4-\{x_1^2+y_1^2=C_1, y_2=0\}$ by setting 
\eqs \tilde{J}= 
\begin{cases}
F_*(J_0) & \text{on } \op{Im} F( \{ \| (s,t) \| \leq 3 \tilde{\e} \} ) \\
j & \text{otherwise}.
\end{cases}
\eeqs 
It's straightforward to check using \eqref{equation:j-pullback} that $\tilde{J}$ is well-defined and smooth.

\lem \label{lemma:compatible-ac}
	The almost-complex structure $\tilde{J}$ is compatible with the standard symplectic form $\o$.
\elem

\pf
It is enough to prove that $\tilde{J}$ is compatible with $\o$ at points $F(p) \in \R^4-\{x_1^2+y_1^2=C_1, y_2=0\}$ where $p \in \{ \| (s,t) \| \leq 2 \tilde{\e} \}$.  One first observes that the splitting $T_{\phi(p)} \R^4 = \op{span} \{ \d_{x_1}, \d_{y_1} \} \oplus \op{span} \{ \d_{x_2}, \d_{y_2} \}$ induces a splitting $${(F_*(J_0))}_{\phi(p)}={J_0^1} \oplus {J_0^2}.$$ Observe that $\o$ also splits as $\o= \o^1 \oplus \o^2$. Hence we only need to check that $J_0^k$ is compatible $\o^k$ for $k=1,2$. This is true for dimension reasons (since any almost-complex structure on a $2$-dimensional symplectic vector space whose induced orientation agrees with the orientation induced by the symplectic form is automatically compatible with that symplectic form).
\epf

Let us now switch gears and discuss a general procedure for constructing a canonical almost-complex structure on the cotangent bundle of a Riemannian manifold. This procedure will be useful to us in the next section. It is originally due to Sasaki and we refer the reader to \cite[Sec.\ 1.3--1.4]{paternain} for a detailed exposition.\footnote{Some conventions in \cite{paternain} are different from ours: in particular, the symplectic form on the cotangent bundle and the metric-induced almost-complex structure both differ by a sign. We also note that although \cite{paternain} mainly considers tangent bundles, all of the relevant constructions commute with the musical isomorphisms so this distinction is entirely superficial.}

Let $(M, g)$ be a Riemannian manifold. The Levi-Civita connection induces a splitting $$T T^*M= \mc{H} \oplus \mc{V},$$ where $\mc{H}, \mc{V}$ are respectively the horizontal and vertical distributions. Given $\t \in T^*M$, let us consider a pair of linear maps: 
\begin{align*} (d_{\t} \pi)^{\flat}: T_{\t} T^*M \to T_{\pi(\t)}^*M, &\hspace{1cm}  K_{\t}: T_{\t} T^*M \to T_{\pi(\t)}^*M. \end{align*}

Here $(d_{\t} \pi)^{\flat}$ denotes the differential of the canonical projection $\pi: T^*M \to M$, composed with the musical isomorphism $TM \to T^*M$. The map $K_{\t}$ is the \emph{connection map} and is defined as follows. Given $\theta \in T^*M$ and $\xi \in T_{\t} T^*M$, choose a path $\g: (-1,1) \to T T^* M$ such that $\g(0)= \t$ and $\dot{\g}=\xi$. Letting $\a:= \pi \circ \g$, we can write $\g=(\a(t), Z(t))$, where $Z$ is a covector field along $\a$. Now define $$K_{(q,p)}(\theta):= \nabla_{\dot{\a}} Z(0).$$

\lem The maps $(d_{\t} \pi)^{\flat}$ and $K_{\t}$ are linear and surjective. Moreover, we have $\op{ker} (d_{\t} \pi)^{\flat} = \mc{V}_{\t}$ and $\op{ker} K_{\t}= \mc{H}_{\t}$. \elem 
\qed

Thus we obtain identifications $(d_{\t} \pi)^{\flat}: \mc{H}_{\t} \to T^*_{\pi(\t)} M$ and $K_{\t}: \mc{V}_{\t} \to T^*_{\pi(\t)} M$. Writing $$T_{\t} T^*M \ni \xi = (\xi_1, \xi_2) \in T^*_{\pi(\t)} M \oplus T^*_{\pi(\t)} M,$$ where the identification $\xi = (\xi_1, \xi_2)$ is induced by $((d_{\t} \pi)^{\flat}, K_{\t})$, we define 
\begin{align*} J_g: T T^*M &\to T T^*M \\ 
J_g(\xi_1, \xi_2)&= (\xi_2, - \xi_1).
\end{align*}

This almost-complex structure is compatible with $\o= d\l$. Observe that if $g$ is Euclidean metric on $\R^{n}$, then the above construction just gives back the standard integrable complex structure on $T^* \R^{n}\simeq \R^{2n}$. 

The following lemma will be useful in the next section. 

\lem \label{lemma:convex-hyp} For $r>0$, hypersurface $\mc{S}_r:= \{ x \in T^*M \mid \|x\|_g = r \}$ is pseudoconvex (see \cite[Sec.\ 2.3]{cieliebak-eliashberg}) for the almost-complex structure $J_g$. \elem
\pf 
Consider the function $H: T^*M \to \R$ defined by $H(q, p):= \frac{1}{2} \langle p, p \rangle_g$. It is an easy exercise (see \cite[Sec.\ 2.8]{cieliebak-eliashberg}) to show that the condition for the level sets of $H$ to be pseudoconvex is the same as the condition for $\op{ker} dH \circ J_g $ to be a contact structure. But note that given $\xi \in T_{(q,p)} T^*M$, we have $dH\circ J_g(\xi)= \langle K(J \xi), p \rangle_g= \langle (d_{\t} \pi)^{\flat} (\xi), p \rangle_g = \l_{can}(\xi)$; cf.\ \cite[Prop.\ 1.21 and Def.\ 1.23]{paternain}. Hence $\l_{can}= dH\circ J_g$, so it suffices to check that $\l_{can}$ restricts to a contact structure on $\mc{S}_r$. This is in turn a consequence of the fact that the radial Liouville vector field is transverse to $\mc{S}_r$. 
\epf

\section{Proof of \Cref{theorem:main}} \label{section:proof-main-thm}

As in the statement of \Cref{theorem:main}, we let $\S$ be a non-compact Riemann surface of finite type and genus $g \geq 0$ (i.e. $\S$ is obtained by removing a positive and finite number of points from a genus $g$ surface). We let $L \sub (T^* \S, d\l_{can})$ be a Lagrangian submanifold which is diffeomorphic to $\R^2$ and agrees outside a compact set with a fiber $F_x$ for some $x \in \S$.

We let $(\R^4, \o)$ be the standard symplectic vector space, where $\o= dx_1 \wedge dy_1 + dx_2 \wedge dy_2$ with respect to the coordinates $(x_1, y_1, x_2, y_2)$. 

\rmk[Notation] 
If $A$ and $B$ are sets, we define $A-B:= \{x \in A \mid x \notin B\}$. In particular, this notation is not reserved for the special case where $B$ is a subset of $A$.
\ermk
 
\subsection{Guide to the proof}\label{subsection:guide-to-the-proof}
The proof of \Cref{theorem:main} is divided into four sections. To help orient the reader, we briefly explain how these sections fit together. The reader may also wish to consult \Cref{subsection:remarks-proof-intro} where the general strategy of the proof of \Cref{theorem:main} is summarized at a more conceptual level. 

In \Cref{section:standard-model} we begin the setup of the proof. In particular, we construct a convenient model for $\S$. We also construct a special Riemannian metric $g_\eta$ on $\S$. We remark that the constructions in \Cref{section:standard-model} do not involve the Lagrangian $L$ (however, certain constants appearing in this section are fixed later in way which does depend on $L$).

In \Cref{subsection:geometric-setup}, we continue the setup of the proof. In particular, we introduce a hypersurface ``with holes" called $W'$, which contains $L$ and is diffeomorphic to $\R^3$ with two solid cylinders removed. We also construct a certain almost-complex structure $J$ (defined on $T^*\S$ with two empty cylinders removed) which depends on the metric $g_\eta$ (and on some other choices). 

\Cref{subsection:filling-planes} is the heart of the proof. Here we construct a symplectic hypersurface $Q$ diffeomorphic to $\R^3$ which contains $L$. This hypersurface is constructed by ``filling in" the holes of $W'$. To do this, we study certain moduli spaces of $J$-holomorphic planes (where $J$ is the almost-complex structure constructed in \Cref{subsection:geometric-setup}). We show that these moduli spaces can be compactified to form embedded solid cylinders which can be glued smoothly to $W'$.  The key compactness statement is \Cref{lemma:closed-image}, which depends heavily on the way in which we have constructed $J$ (which, in turn, depends on $g_\eta$). It is necessary to use both $J$-convexity and positivity of intersection to confine the curves to a priori given subsets of $T^*\Sigma$, after which we are able to allude to the SFT-compactness theorem. The gluing happens in \Cref{proposition:smoothing-procedure}, which is similar to \cite[Sec.\ 5.3]{dgi} and \cite[Lem.\ 5.5]{wendlduke}. 

In \Cref{subsection:completion-proof}, we complete the proof of the theorem. We construct a hypersurface $\tilde{Q}$ as a certain deformation of $Q$. As explained in \Cref{subsection:remarks-proof-intro}, we then use $\tilde{Q}$ to construct the desired Hamiltonian isotopy from $L$ to a cotangent fiber. \Cref{subsection:completion-proof} does not involve any pseudoholomorphic curves, and the arguments are similar to those in \cite{eli-pol2}. 

We note that the analytical arguments (e.g.\ ``filling in" the holes of a hypersurface with holomorphic planes) in the proof of \Cref{theorem:main} are quite general, and could plausibly be used to study related problems on other symplectic $4$-manifolds of interest.

\subsection{A standard model for $\S$} \label{section:standard-model}
 
It will be convenient to work with a standard model for the abstract Riemann surface $\S$. Since $\S$ is of finite type, its Euler characteristic is well-defined and given by the formula $\chi(\S)= 2 -2g-p$, where $p \geq 1$ is the number of punctures. 

For $j=1,2,..., 2(2g+p-1)$, let $$I_j = \{ (y_1, y_2) \in \R^2 \mid y_1 \in [j-1/4, j+1/4], y_2=-1 \}.$$ 

For $k=1,2,\dots, 2g+p-1$, let $S_k$ be an abstract manifold with corners equipped with an identification 
\eq \label{equation:sk-identification} S_k= \{(u_1, u_2) \in \R^2 \mid |u_1| \leq 1/4, |u_2| \leq 1\}. \eeq Let $S_k^{\pm} = \{|u_1| \leq 1/4, u_2=\pm 1\} \sub S_k$. We view $S_k^{\pm}$ as oriented $1$-manifolds whose orientation is inherited from the standard orientation on $S_k$. 

Let $\S_+= \{(y_1, y_2) \in \R^2 \mid y_2 > -1 \}$ and let $$\psi: \bigsqcup_{k=1}^{2g+p-1} (S_k^+ \sqcup S_k^-) \to  \bigsqcup_{j=1}^{2(2g+p-1)} I_j \sub \ov{\S}_+$$ be an orientation-preserving diffeomorphism.  We now write 
\eq \label{equation:standard-model} \ov{\S}= \ov{\S}_+ \bigcup_{\psi} \lt( \bigsqcup_{k=1}^g S_k  \rt) \eeq
and define 
\eq \label{equation:sigma} \S= \ov{\S}^{\circ}. \eeq

It will be convenient to assume that $\psi$ extends near each $S_k^{\pm}$ to a Euclidean isometry (i.e. translation and rotation) with respect to the coordinates \eqref{equation:sk-identification}. Hence the standard Euclidean metric descends to $\S$ under the gluing map $\psi$. We let $g_e^{\S}$ denote this metric. 

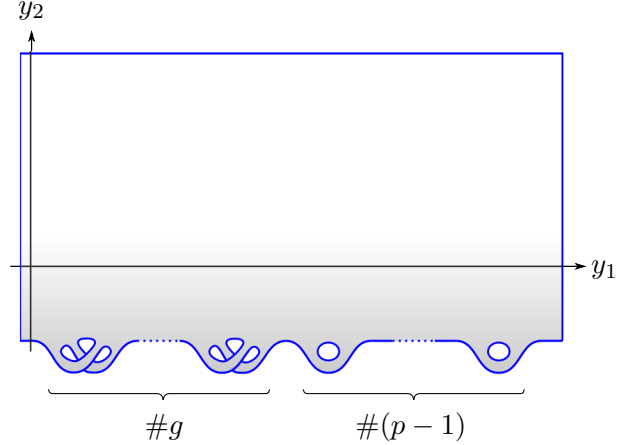
\begin{figure}
\centering
\vspace{3mm}
\labellist
\pinlabel $y_1$ at 224 49
\pinlabel $y_2$ at 8 147
\pinlabel $\#g$ at 58 -10
\pinlabel $\#(p-1)$ at 151 -10
\endlabellist
\includegraphics{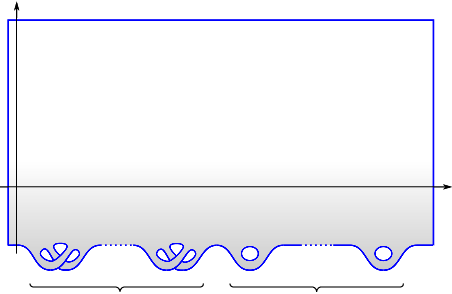}
\vspace{5mm}
\caption{Any compact connected surface with nonempty boundary can be realized as shown in the figure (here the boundary has smoothable corners), where $g \ge 0$ is the genus and $p \ge 1$ is the number of components of its boundary}
\label{fig:surface}
\end{figure}

It is elementary to show that every non-compact Riemann surface of finite-type is diffeomorphic to the Riemann surface which results from the above construction, for an appropriate choice of $\psi$; see \Cref{fig:surface}. Hence there is no loss of generality in taking \eqref{equation:standard-model} and \eqref{equation:sigma} as a model for $\S$. For the remainder of this section, we therefore assume that $\psi$ is fixed and that $\S$ is defined by \eqref{equation:standard-model} and \eqref{equation:sigma}. We may moreover assume without loss of generality that $L$ agrees outside a compact set with the fiber $F_0$ over the point $(0,0) \in \S_+ \sub \S$. 

Observe that there is an identification of symplectic manifolds 
\begin{align} \label{equation:id1} 
(T^* \S_+, d\l_{can}) &\simeq \{(x_1, y_1, x_2, y_2) \mid y_2 > -1 \} \sub (\R^4, \o=dx_1\wedge dy_1+ dx_2 \wedge dy_2) \\
(b_1 dy_1 + b_2 dy_2)_{(a_1, a_2)} &\mapsto (b_1, a_1, b_2, a_2). \nonumber
\end{align}
 
We will routinely make this identification in the sequel without further note. 

Finally, we let \eq \label{equation:boundary} B:= \ov{\S}- \S. \eeq  Note that $B$ is naturally a piecewise smooth $1$-manifold, since $\ov{\S}$ is a manifold with corners. We let $B_0, B_1, \dots, B_{\beta}$ be an enumeration of the components of $B$ for $\beta \geq 0$. After possibly relabeling, we may assume that $B_0$ is homeomorphic to $\R$ while the other components are homeomorphic to $S^1$. 

Given $\eta< 1/100$ and $1 \leq k \leq \beta$, let $$\phi^{\eta}_k: \R/\Z \tms (0, \eta) \to \S$$ be a smooth embedding which extends to a continuous embedding $(\R/\Z \tms [0,\eta), \R / \Z \tms \{0 \}) \to (\S \cup B_k, B_k)$. Let $\phi^{\eta}_0: \R \tms (0, \eta) \to \S$ be a smooth embedding which extends to a continuous embedding $(\R \tms [0, \eta), \R \tms \{0 \}) \to (\S \cup B_0, B_0)$, and such that 
\eq\label{equation:euclidean-large-scale} \phi^{\eta}_0(x,t) = (x, t-1) \eeq 
for $|x|$ large enough. The maps $\phi^{\eta}_k$ define a ``collar" around the boundary component $B_k$.  We let 
\eq \label{equation:eta-collar} \S_{\eta}= \bigcup_{k=0}^{\beta} \op{Im} \phi^{\eta}_{k} \eeq be the union of the collars. 

Finally, let us fix Riemannian metrics $g_k$ (for $0 \leq k \leq \beta$) on $\op{im} \phi^{\eta}_k$ having the following properties:
\begin{itemize}
\item[(i)] For $k \geq 1$ (resp. $k=0$), we have $g_k=g_{e}^{\S}$ on the set $\phi^{\eta}_k(\R /\Z \tms (3\eta/4, \eta))$ (resp. on $\phi^{\eta}_0( \R \tms (3\eta/4, \eta))$), where $g_{e}^{\S}$ is the Euclidean metric on $\S$ which was introduced above. (Recall that $g_{e}^{\S}$ is well-defined on $\S$ in view of our choice of gluing map $\psi$).
\item[(ii)] For $k \geq 1$ (resp. $k=0$), we have $g_k=(\phi^{\eta}_k)_*g_{e}$ on the set $\phi^{\eta}_k(\R /\Z \tms (0, \eta/4))$ (resp. on $\phi^{\eta}_0( \R \tms (0, \eta/4))$), where $g_{e}$ denotes the standard Euclidean metric on $ \R /\Z \tms (0, \eta/4)$ (resp. on $ \R \tms (0, \eta/4)$). 
\item[(iii)] For $|x|$ large enough, $g_0=(\phi^{\eta}_0)_*g_{e}$.%, i.e.\ on the set $\{(x,y) = \phi^{\eta}_0((x,y))$ for $|x|$ large enough. 
\end{itemize}

We let $g_{\eta}$ be a Riemannian metric on $\S$ defined by setting 

\begin{align}\label{align:glued-metric}
g_{\eta}= \begin{cases} g_k & \text{on } \S_{\eta}, \\
g_{e}^{\S} & \text{on } \S - \S_{\eta}.
\end{cases}
\end{align}

Condition (i) above ensures that $g_{\eta}$ indeed defines a smooth metric on $\S$. Condition (ii) says that $g_\eta$ has a nice ``product" structure near the ends of $\S$. Condition (iii), along with \eqref{equation:euclidean-large-scale}, implies that $g_\eta= g_e^{\S}$ outside a compact set. 

The metric $g_\eta$ will be important in the next sections. As explained in \Cref{subsection:guide-to-the-proof}, we are going to construct an almost-complex structure $J$ on $T^*\S$ minus two Lagrangian cylinders  $\Phi(\mc{L}^+ \cup \mc{L}^-)$ (where $\Phi, \mc{L}^\pm$ will be defined shortly). Away from the Lagrangian cylinders, $J$ will look like the Sasakian almost-complex structure induced by $g_\eta$ (in contrast, near the Lagrangian cylinders, $J$ will be cylindrical). The properties of $g_\eta$ will be important when we study $J$-holomorphic curves on $T^*\S- \Phi(\mc{L}^+ \cup \mc{L}^-)$. In particular, the fact that $g_\eta$ has a nice flat product structure (i.e.\ condition (ii)) in a neighborhood of $\ov{\S} - \S$ will be crucial for ensuring that $J$-holomorphic curves in $T^*\Sigma$ cannot approach the fibres above $\overline{\Sigma} - \Sigma$. Condition (iii) will similarly be used to prevent $J$-holomorphic curves from escaping in the direction $|y_1| \to \infty$. (See the proof of \Cref{lemma:closed-image} for details.)

\subsection{Construction of a hypersurface ``with holes"} \label{subsection:geometric-setup}

In this section, we continue setting up the proof of \Cref{theorem:main}. In particular, we construct a ``hypersurface with holes" \eqref{equation:w'} and an almost-complex structure $J$ on $T^*\S-\Phi(\mc{L}^+ \cup \mc{L}^-)$ (where $\Phi(\mc{L}^+ \cup \mc{L}^-)$ are Lagrangian cylinders to be defined shortly). Both of these constructions will play an important role in \Cref{subsection:filling-planes}.

We begin with the construction of our ``hypersurface with holes". Recall from the previous section that $F_0 \sub T^*\S$ denotes the fiber over the point $(0,0) \in \S_+ \sub \S$ via the identification \eqref{equation:id1}. Given $\e$ small enough, the symplectic neighborhood theorem provides a symplectic embedding $$\Phi_0: \op{Op}_{\e}(F_0) \to T^*\S$$ with the property that:
\begin{itemize}
\item $\Phi_0(F_0)= L$,
\item $\Phi_0$ restricts to the identity on $(T^* \S - \mc{K}) \cap \op{Op}_{\e}(F_0)$, where $\mc{K} \subset T^*\Sigma$ is a \emph{compact} connected embedded closed submanifold of codimension zero with smooth boundary.
\end{itemize}
After possibly making $\mc{K}$ larger, we may assume that the image of $\mc{K} \cap \op{Op}_{\e}(F_0)$ under $\Phi_0$ is entirely contained in $\mc{K}$. After possibly making $\eta>0$ smaller, we may in addition assume that $\mc{K} \cap \pi^{-1}(\S_{\eta}) = \emptyset$, where $\pi: T^* \S \to \S$ is the canonical projection; see \eqref{equation:eta-collar} for the definition of $\S_{\eta}$. 

Let us fix $M >0$ large enough so that the inclusion
\begin{equation}\label{equation:M-in-base}
T^*\S_+ \cap \mathcal{K} \subset \{|x_i| <M/2, |y_i| < M/2, i=1,2\} \sub T^*\S_+
\end{equation}
is satisfied. (\eqref{equation:M-in-base} will be used later in the proof of \Cref{lemma:closed-image}.)

We view $\e, \Phi_0, \mc{K}, \eta, M$ as fixed for the remainder of this section. 

\begin{figure}
\begin{tikzpicture}[x=0.6cm,y=0.6cm]
\clip(-8.5,-8.5) rectangle (8.5,8.5);
\draw [line width=2pt,domain=-8.5:8.5] plot({(-0-0*\x)/3},\x);
\draw [line width=1.2pt] (-8.5,0) -- (8.5,0);
\draw [line width=1.2pt] (25.800166666666655,0) circle (25.370166666666652);
\draw [line width=1.2pt] (-25.800166666666655,0) circle (25.370166666666652);
\draw [->,line width=1.2pt] (0,0) -- (8.5,0);
\draw [->,line width=1.2pt] (0,0) -- (0,8.5);
\draw (8,-0.5) node {$y_1$};
\draw (-0.5,8.0) node {$x_1$};
\draw (0.5,-7.0) node {$F_0$};
\draw (1.5,-5.0) node {$\mathcal{L}^+$};
\draw (-1.5,-5.0) node {$\mathcal{L}^-$};
\end{tikzpicture}
\caption{A slice $\{(x_1,y_1,x_2, y_2 \mid x_2=\text{constant}, y_2= 0\}$.}
\label{figure-circles}
\end{figure}
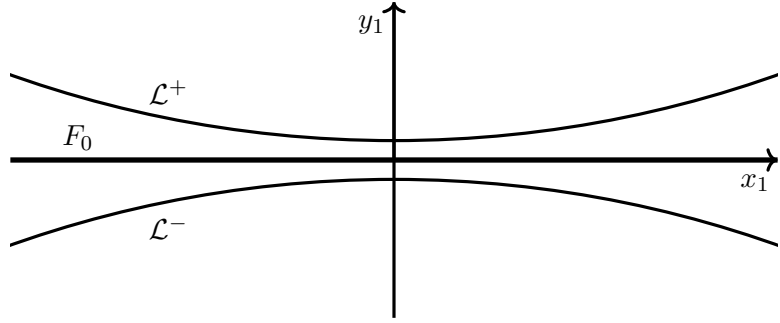

Define the  open set 
\eq X:= \lt(T^* \S-\mc{K} \rt) \bigcup  \op{Op}_{\e}(F_0), \eeq 
and let $\Phi: X \to T^*\S$ be defined by 
\eqs
\Phi(p)= \begin{cases} \Phi_0(p) & \text{if } p \in \op{Op}_{\e}(F_0) \\
p & \text{otherwise}.
\end{cases}
\eeqs
One readily checks that $\Phi$ is well-defined. (Indeed, this follows from the fact that $\Phi_0 (\mc{K}\cap \op{Op}_\e(F_0))$ is contained in $\mc{K}$, and that $\mc{K} \cap \pi^{-1}(\Sigma_\eta)= \emptyset$.) Note also that 
\eq \pi^{-1}(\S_{\eta}) \sub X.\eeq
Choose $C$ large enough (depending on $\e, M$) so that the sets
\eq\label{equation:wl} 
W_{l}= \{ x_1^2+ (y_1+ C)^2 \geq (C-\e/2)^2, x_1^2+ (y_1- C)^2 \geq (C-\e/2)^2,  x_2=l, y_2= 0 \}
\eeq
are contained in $X$ for all $l \in \R$. In particular, this means that the Lagrangian cylinders
$$
\mc{L}^{\pm}:= \{ x_1^2 + (y_1 \mp C)^2 = (C- \e/2)^2, y_2=0 \}
$$ 
are also contained in $X$. In fact, as illustrated in \Cref{figure-circles}, the subset $\mc{L}^{\pm} \cap \mc{K}$ of the Lagrangian $\mc{L}^\pm$ is contained inside an $\epsilon$-small Lagrangian push-off of the fibre $F_0$, and is hence contained inside the subset $\op{Op}_{\e}(F_0) \subset X$.

For future reference, let us define 
\eq\label{equation:def:w'l} W'_l:= \Phi(W_l)\eeq 
for $l \in \R$. We now introduce our promised \emph{hypersurface with holes:} 
\eq\label{equation:w'} W':= \bigcup_{l \in \R} W'_l.\eeq
In conclusion: $W'$ is a symplectic hypersurface which contains $L$ and is diffeomorphic to $\R^3$ with two solid cylinders drilled out. The boundary of $W'$ is a pair of Lagrangian cylinders $\Phi(\mc{L}^+ \cup \mc{L}^-)$. In the next section, we will glue in two solid cylinders constructed from moduli spaces of pseudoholomorphic planes. 

We now turn to the construction of the almost-complex structure $J$ on $T^*\S- \Phi(\mc{L}^+ \cup \mc{L}^-)$. To begin with, define $\tau^{\pm}: \R^4 \to \R^4$ by $\tau^{\pm}_i(x_1,y_1,x_2, y_2)= (x_1,y_1\pm  C,x_2, y_2)$. Let $C_1:= C-\e/2$ and define 
$$F^{\pm}:= \tau^{\pm} \circ F.$$ We remind the reader that our definition of $F$, which was stated in \eqref{align:def-of-F}, depends on $C_1$. It follows by definition that $F^{\pm}$ takes the zero section $0_{S^1 \tms \R} \sub T^*(S^1_{\theta} \tms \R_r)$ to the Lagrangian cylinder $\mc{L}^\pm$. Thus $\Phi \circ F^\pm$ defines coordinates on a Weinstein neighborhood of $\Phi(\mc{L}^\pm)$.

Let us now return to the metric $g_{\eta}$ constructed at the end of \Cref{subsection:a-c-structures}. Noting that $\eta$ is now fixed, we let $j_S$ be the ``Sasakian" almost-complex structure on $T^*\S$ induced by the metric $g_{\eta}$, as described at the end of \Cref{subsection:a-c-structures}. Observe that $j_S$ restricts to the standard integrable almost-complex structure on $\R^4$ away from $\pi^{-1} (\S_{\eta})$, since $g_{\eta}$ is just the standard Euclidean metric on $\S- \S_{\eta}$ (here, we are as usual using the identification \eqref{equation:id1}). 

Next, fix $\tilde{\e}$ small enough (depending on $\e, M, C$) so that $$F^{\pm}( \{ \| (t,s ) \| \leq 3 \tilde{\e} \} ) \sub X,$$ and define an almost-complex structure $J'$ on $T^*\S- (\mc{L}^+ \cup \mc{L}^-)$ by setting
\eq
J' = \begin{cases}
({\tau^{\pm}})_* \tilde{J} & \text{on } F^{\pm}_i( \{0 < \| (t,s ) \| \leq 3 \tilde{\e} \} ) \\
j_{S} & \text{otherwise}.
\end{cases}
\eeq
Using \Cref{lemma:compatible-ac}, it is straightforward to verify that $J'$ is well-defined and compatible with $\o$. 

Finally, we fix an almost-complex structure $J$ on $T^*\S - \Phi(\mc{L}^+ \cup \mc{L}^-)$ which is compatible with $d \l_{can}$, and such that $J= \Phi_* J'$ in $\Phi(X- (\mc{L}^{+} \cup \mc{L}^-))$. 

It follows from \Cref{lemma:cylindrical} and the definition of $J$ that this almost-complex manifold has negative cylindrical ends around the $\Phi(\mc{L}^{\pm})$ of the form 
\eq\label{equation:cylindrical-ac-ref} ((-\infty, 0] \tms S_{\tilde{\e}}, J_{cyl}),\eeq 
where $J_{cyl}$ was defined in the paragraph following \Cref{lemma:cylindrical}. 

In conclusion: our almost-complex structure $J$ is cylindrical near $\Phi(\mc{L}^+ \cup \mc{L}^-)$, and it looks like the Sasakian almost-complex structure $j_S$ induced by $g_{\eta}$ away from $\Phi(\mc{L}^+ \cup \mc{L}^-) \cup \mc{K}$. Inside $\mc{K}$ (which is compact), $J$ is a mystery almost-complex structure which is compatible with $d \l_{can}$.

We have now finished assembling the ingredients for setting up the proof of \Cref{theorem:main}.

\subsection{Construction of a symplectic hypersurface} \label{subsection:filling-planes}

In this section, we will fill-in the holes of $W'$ by gluing certain (appropriately compactified) moduli spaces of $J$-holomorphic planes. More precisely, we will consider punctured holomorphic curves in the almost-complex manifold $$(T^*\S - \Phi(\mc{L}^+ \cup \mc{L}^-), J).$$ 

For $\s \in ((-\infty, -M] \cup [M, \infty))$, let 
\eq\label{equation:standar-linear-plane} u^{\pm}_{0, \s}: \C \to T^*\S - \Phi(\mc{L}^+ \cup \mc{L}^-)\eeq be a $J$-holomorphic plane whose image is the set $\{ x_1^2 + (y_1\mp C)^2 < (C-\e/2)^2, x_2=\s, y_2= 0 \}$. Such a plane is unique up to reparametrization. 

Let $\mc{M}^{\pm}$ be the connected component of the moduli space of unparametrized $J$-holomorphic planes in $(T^*\S - \Phi(\mc{L}^+ \cup \mc{L}^-))$ containing $u^{\pm}_{0,M}$. We compute the index of this component. 

\lem \label{lemma:index-computation}
Given a plane $u \in \mc{M}^{\pm}$, we have $\op{ind}(u)= 1$ and $\op{ind}^{c}(u)=0$, where $\op{ind}^c(-)$ is the index for the moduli problem for planes constrained to be asymptotic to a fixed Reeb orbit (see \Cref{subsection:punctured-curves}).
\elem

\pf
First of all, if $u$ is a pseudoholomorphic plane with asymptotic orbits in a Morse--Bott submanifold $X$, then it follows from \eqref{equation:index} and \eqref{equation:index-difference} that $\op{ind}(u)- (\op{dim}(X)-1)= \op{ind}^{c}(u)$. Hence it is enough for us to prove that $\op{ind}(u)=1$. Since the index is invariant under homotopies, we only need to compute the index of $u_{0,M}^{\pm}$. 

To compute $\op{ind}(u)=1$ one could use the formula $\op{ind}(u)=-1+\mu_{\mc{L}^\pm}[\overline{u}]$ (see e.g.~\cite[Sec.\ 3.1]{dgi}), where $\mu_{\mc{L}^\pm}[\overline{u}]$ denotes the Maslov index of the relative cycle $[\overline{u}] \in H_2(T^*\Sigma,\mc{L}^\pm)$ obtained by compactification.

Alternatively, one can perform the following direct computation. By combining \Cref{definition:index} and \Cref{lemma:cz-index}, one finds that it is enough to verify that $c_1^{\Theta}(u_{0,M}^{\pm})=1$, where $\Theta=\{ -J_{cyl}\d_\phi\}$ is the complex trivialization of the contact planes defined in \Cref{lemma:cz-index} (the coordinate $\phi$ is the angular coordinate on the unit cotangent fibre; see \Cref{subsection:a-c-structures}). Note that the vectors $\d_\phi$ tangent to the contact distribution on $S_{\tilde{\e}}$ are identified with vectors inside $\R_{>0} \cdot  \partial_{y_2}$ when restricted to the subset $\{y_2=0\}$ of either cylindrical end; see the definition of the maps $F^\pm$ above, together with the definition of $F$ in \Cref{subsection:a-c-structures}.  
	
We are now ready to compute $c_1^{\Theta}(u_{0,M}^{\pm})=1$. Consider the parametrization of the pseudoholomorphic plane
\begin{gather*}
u_{0,M}^{\pm} \colon B^2(\sqrt{C -\epsilon/2}) \to T^*\Sigma \setminus \{\mc{L}^+ \cup \mc{L}^-\},\\
(r^{\pm}, \t^{\pm}) \mapsto (r^\pm, \t^\pm, x_2=M, y_2=0),
\end{gather*}
defined by polar coordinates $r^{\pm}=(x_1^2+(y_1 \mp C)^2)^{1/2}$ and $\t^{\pm}= \tan^{-1}((y_1 \mp C)/ x_1)$; in other words $(r^\pm,\t^{\pm})$ are polar coordinates on the $(x_1,y_1)$-plane centred at the point $\{(x_1,y_1)=(0,\pm C)\}$. Then $\d_{r^{\pm}} \wedge \d_{y_2}$ is a section of the complex determinant bundle $\Lambda^2_{\C}((u^{\pm}_{0,M})^*T(T^*\S))$ which is asymptotically constant (up to homotopy) with respect to $\d_{r^{\pm}} \wedge\Theta$, and which vanishes to order $1$ as $r^{\pm} \to 0$. (Note that, along the asymptotic Reeb orbit $r \to \sqrt{C -\epsilon/2}$ of the plane, the vector field $\partial_{y_2}$ defines a complex trivialisation of the contact planes in the same homotopy class as $\Theta$ from \Cref{lemma:cz-index}.) 
\epf

Let $(\mc{M}^{\pm})_{reg} \sub \mc{M}^{\pm}$ be the open subset of transversally cut-out planes. Since the planes under consideration are asymptotic to a primitive closed geodesic and are therefore simply covered, the reparametrization group acts freely. It follows that $(\mc{M}^{\pm})_{reg}$ is a smooth $1$-dimensional (Hausdorff) manifold; see \cite[Thm.\ 0]{wendlautomatic}.  

We let $\mc{U}^{\pm} \sub (\mc{M}^{\pm})_{reg}$ be the connected component containing $u^{\pm}_{0, M}$. 

It will be useful to record the following lemma, which shows that the elements of $\mc{U}^{\pm}$ are also transverse for a constrained moduli problem. 

\lem \label{lemma:restricted-moduli}
Given a plane $u \in \mc{U}^{\pm}$, the linearization of $\overline{\d}_J$ at $u$ through planes whose asymptotic orbit is fixed is surjective.
\elem
\pf
We first view $u$ as a plane with an unconstrained puncture. It's clear that the standard plane $u^{\pm}_{0, M}$ has vanishing Siefring self-intersection number. Indeed, for $t>0$, the plane $u^{\pm}_{0,M+t}$ is disjoint from $u^{\pm}_{0,M}$ and shares no asymptotic orbit with it, so this follows by combining \eqref{item:homotopy-invariance} and \eqref{item:disjointasymptotics} in \Cref{fact:siefring-number}. Since the Siefring intersection number is constant in families (see \Cref{fact:siefring-number} \eqref{item:homotopy-invariance}), the other planes in $\mc{U}^{\pm}$ also have vanishing Siefring self-intersection number. 
	
Note next that these planes must have vanishing normal Chern number: indeed, since the normal Chern number is invariant in families, it is enough to check that $c_N(u^{\pm}_{0,M})=0$. But by definition, we have $c_N(u^{\pm}_{0,M})= \op{ind}(u^{\pm}_{0,M})-2+\#\G_0= 1-2+1=0,$ where we have used \Cref{lemma:index-computation} and \Cref{lemma:cz-index}. It follows from the adjunction formula (\Cref{fact:adjunction}) that all the planes are all embedded. In particular, $u$ is embedded. (Note that in applying the adjunction formula, we use the fact that the planes under consideration are asymptotic to a primitive geodesic, which means that the term $\op{cov}_\infty(-)$ vanishes). 
	
We now view $u$ as a plane with a constrained puncture. We then have $\op{ind}^c(u)=0$ while \Cref{definition:normal-chern} implies that $2c^c_N(u) = \#\G^c_0- 2 \leq -1$. It follows that $c^c_N(u) < \op{ind}^c(u)$, so the desired conclusion follows by automatic transversality (\Cref{fact:automatic-transversality}). (The notation $(-)^c$ is intended to emphasize that we are working with a constrained puncture.)
\epf

Let $\G^{\pm}$ be the manifold of closed, simple and positively oriented geodesics of the Lagrangian cylinder $\Phi(\mc{L}^{\pm}) \sub T^*\S$. There is an identification $$ \g^{\pm}: \R \xrightarrow{\sim} \G^{\pm}$$ which takes $l \in \R$ to the unique geodesic in $\G^{\pm}$ which passes through $\Phi \circ F^{\pm}(0,0,l,0)$. 

Let ${(\op{ev}_{\d})}^{\pm}: \mc{M}^{\pm} \to \G^{\pm}$ be the asymptotic evaluation map (see \cite[A.2.]{wendlduke}) and define
\begin{align*}
\Psi^{\pm}: \mc{U}^{\pm}_i &\to \R \\
u &\mapsto (\g^{\pm})^{-1} \circ {(\op{ev}_{\d})}^{\pm}.
\end{align*}
Since $\mc{U}^{\pm}$ consists of transversally cut out planes by definition, it follows that ${(\op{ev}_{\d})}^{\pm}$ and hence $\Psi^{\pm}$ are local diffeomorphisms; cf. \cite[Prop.\ 5.11(i)]{dgi}.  

\rmk
The target of the asymptotic evaluation map defined in \cite[A.2.]{wendlduke} is in fact a line bundle $E^{\pm}$ over $\G^{\pm}$. Our asymptotic evaluation map is obtained by composing the map in \cite[A.2.]{wendlduke} with the projection $E^{\pm} \to \G^{\pm}$. To verify that the composition is a local diffeomorphism, one needs to establish as in \cite[Prop.\ 5.11(i)]{dgi} that the elements of $\mc{U}^{\pm} $ are transverse for the moduli problem with constrained orbit as well as for the moduli problem where the orbits allowed to vary in the Morse--Bott family. This is why we need \Cref{lemma:restricted-moduli}. 
\ermk

We wish to show that $\Psi^{\pm}$ is in fact a diffeomorphism. This will be deduced from the following sequence of lemmas.

\lem \label{lemma:injective}
The map $\Psi^{\pm}$ is injective.
\elem
\pf
Suppose that $u, v \in \mc{U}^{\pm}$ are asymptotic to the same Reeb orbit $\g$. We wish to show that $u=v$ (up to reparametrization). Suppose for contradiction that $u \neq v$. As observed in the proof of \Cref{lemma:restricted-moduli}, the Siefring self-intersection number of the standard plane $u^{\pm}_{0,M}$ vanishes. Since $\mc{U}^{\pm}$ is connected and contains $u^{\pm}_{0,M}$, it follows that $u*v=0$. On the other hand, an argument originally due to Hind and Lisi (see \cite[Thm.\ 4.2]{hindandlisi} as well as \cite[Lem.\ 5.13]{dgi}) shows that $u* v>0$ if both planes are asymptotic to the same Reeb orbit. The idea is that $u, v$ differ asymptotically by an eigenfunction of the asymptotic operator $\bf{A}_\g$ with positive eigenvalue (this is a consequence of the asymptotic convergence formulas alluded to in \Cref{subsection:punctured-curves}). Such eigenfunctions turn out to have positive winding number, which means that a small pushoff of $v$ will intersect $u$ positively. This proves the claim. \epf

By combining \Cref{lemma:injective} with the previously observed fact that $\Psi^{\pm}$ is a local diffeomorphism, we conclude that $\Psi^{\pm}: \mc{U}^{\pm} \to \R$ is a smooth embedding. 

The same argument as in \Cref{lemma:injective} also gives the following statement, which will be useful later. 
\lem \label{lemma:unique-standard}
Given $| l | \geq M$, the plane $u^{\pm}_{0, l}$ is the unique $J$-holomorphic plane asymptotic to the geodesic $\g^{\pm}(l)$. 
\qed
\elem

Our next task is to show that $\Psi^{\pm}$ has closed image. To this end, it will be useful to introduce certain barriers which provide restrictions on holomorphic planes escaping to infinity. 

\begin{figure} 
\begin{tikzpicture}[x=0.6cm,y=0.6cm]
\clip(-8.5,-8.5) rectangle (8.5,8.5);
\draw [line width=1.2pt,domain=-8.5:8.5] plot({(-0-0*\x)/3},\x);
\draw [line width=1.2pt] (-8.5,0) -- (8.5,0);
\draw [line width=2pt]   (1,8.0) .. controls (1,8) and (1,8) .. (1,5) .. controls (1,5) and (1,3) .. (0,3);
\draw [line width=2pt]   (-1,8) .. controls (-1,8) and (-1,8) .. (-1,5) .. controls (-1,5) and (-1,3) .. (0,3);
\draw [line width=2pt]   (1,-8.0) .. controls (1,-8) and (1,-8) .. (1,-5) .. controls (1,-5) and (1,-3) .. (0,-3);
\draw [line width=2pt]   (-1,-8.0) .. controls (-1,-8) and (-1,-8) .. (-1,-5) .. controls (-1,-5) and (-1,-3) .. (0,-3);
\draw [->,line width=1.2pt] (0,0) -- (8.5,0);
\draw [->,line width=1.2pt] (0,0) -- (0,8.5);
\draw (8,-0.5) node {$y_2$};
\draw (-0.5,8.0) node {$x_2$};
%\draw (-7.0,0.5) node {$F_0$};
\draw (1.25,-3.25) node {$\mathcal{H}^-$};
\draw (1.25,3.25) node {$\mathcal{H}^+$};
\draw (1,1) node {$\mathcal{V}$};
\draw (0.5,-5.0) node {$\mathcal{V}^-$};
\draw (0.5,5.0) node {$\mathcal{V}^+$};
\end{tikzpicture}
\caption{The projection of the Levi-flat hypersurface $\mc{H}$ to the $(x_2,y_2)$-plane.}
\label{figure:hypersurface}
\end{figure}
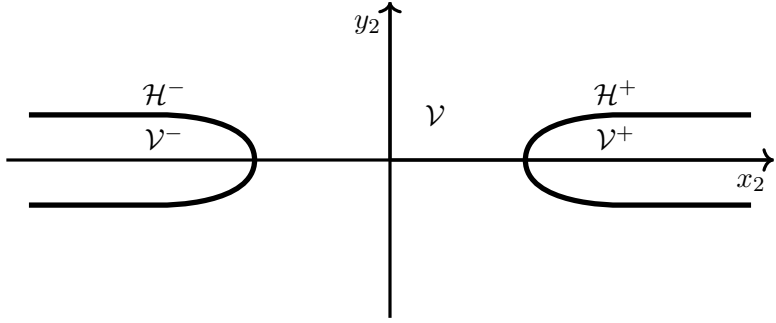

Let $f: \R_{\geq 0} \to [0,1/2] \sub \R$ be a non-decreasing function with the property that $f(x) = \sqrt{x}$ for $x \in [0, 1/2]$. Let us now consider the hypersurfaces (see \Cref{figure:hypersurface})
\begin{align*}
\mc{H}^{+} &= \{ (x_1, y_1, x_2, y_2) \sub T^*\S_+ \mid x_2  \geq 2M, |y_2| = f(x_2-2M) \}, \\
\mc{H}^{-} &= \{ (x_1, y_1, x_2, y_2) \sub T^*\S_+ \mid x_2  \leq -2M, |y_2| = f(2M-x_2) \}.
\end{align*}

Observe that we can view $\mc{H}^{\pm}$ as being contained in the set $\{(x_1, y_1, x_2, y_2) \mid y_2 \geq -1/2\} \sub T^*\S_+ \sub T^*\S$.  We let $\mc{H}= (\mc{H}^+ \cup \mc{H}^-) - \Phi(\mc{L}^+ \cup \mc{L}^-)=  (\mc{H}^+ \cup \mc{H}^-) - (\mc{L}^+ \cup \mc{L}^-)$. 
Let us record the following important lemma.

\lem \label{lemma:hypersurface-separating}
The hypersurface $\mc{H} \subset T^*\S - \Phi(\mc{L}^+ \cup \mc{L}^-)$ is properly embedded and foliated by properly embedded, boundaryless $J$-holomorphic surfaces (in particular it is Levi-flat). Moreover, $T^*\S - (\mc{H}^+ \cup \mc{H}^-)$ consists of three connected components, and hence so does $(T^*\S - \Phi(\mc{L}^+ \cup \mc{L}^-)) - \mc{H}$. We can uniquely label these components $\mc{V}^+, \mc{V}, \mc{V}^-$ by requiring that $\mc{H}^{+}$ (resp. $\mc{H}^-$) separates $\mc{V}$ and $\mc{V}^+$ (resp.\ $\mc{V}$ and $\mc{V}^-)$. 
\elem
\pf
Observe that the planes $\{ (x_1, y_1, x_2, y_2) \mid y_2= \pm f(x_2-M) \}$ are $J$-holomorphic, for each fixed $(x_2, y_2)$. This is immediate from the observation (cf.\ \Cref{lemma:compatible-ac}) that $J$ preserves the splitting $\{\d_{x_1}, \d_{y_1} \} \oplus \{ \d_{x_2}, \d_{y_2} \}$ in the region $\{(x_1, y_1, x_2, y_2) \mid y_2 \geq -1/2, |x_2| \geq M \} \sub T^*\S_+ \sub T^*\S$. 
The fact that $T^*\S - (\mc{H}^+ \cup \mc{H}^-)$ consists of three connected components is clear from the definition of $\mc{H}^+ \cup \mc{H}^-$. 
\epf

We now proceed to the key compactness result of this section. As a small piece of notation, given $s \in \op{Im} \Psi^{\pm} \sub \R$, let us write $u^{\pm}_s:= {(\Psi^{\pm})}^{-1}(s)$.

\lem \label{lemma:closed-image}
The image of $\mc{U}^{\pm}$ under $\Psi^{\pm}$ is closed as a subset of $\R$.
\elem

Before proving this lemma, let us first discuss informally what could go wrong. We are going to consider a sequence of $J$-holomorphic planes in $T^*\S- \Phi(\mc{L}^+ \cup \mc{L}^-)$, and we want to prevent this sequence from escaping to infinity. A priori, there are multiple different ways in which these planes could escape: (i) planes could escape along the ends of $\S$ (which is of course a \emph{non-compact} Riemann surface); (ii) planes could escape along the fibers of $\S$; (iii) planes could escape into $\Phi(\mc{L}^+ \cup \mc{L}^-)$. Ultimately, we want to argue that the image of our planes under the inclusion $T^*\S- \Phi(\mc{L}^+ \cup \mc{L}^-) \hookrightarrow T^*\S$ remains in a compact subset. Then, the conclusion of the lemma will follow rather directly using the SFT compactness theorem.

To prevent planes from escaping along the ends of $\S$, we will use our construction of $g_\eta$ in \Cref{section:standard-model}. More precisely, near the ends of $\S$, $J$ looks like the Sasakian metric induced by $g_\eta$. By construction, $g_\eta$ has a flat product structure near these ends. As a result, the ends are foliated by $J$-holomorphic curves, which in turn allows us to conclude (by a positivity of intersection argument) that our sequence of planes cannot escape. 

To prevent planes from escaping along the direction of the fibers of $\S$, the argument is more delicate. This is because the almost-complex structure is not Sasakian everywhere; indeed, near $\Phi(\mc{L}^+ \cup \mc{L}^-)$, the almost-complex structure is \emph{cylindrical}, whereas far away from $\Phi(\mc{L}^+ \cup \mc{L}^-)$, it becomes Sasakian. We therefore need to argue differently in each case: far away from $\Phi(\mc{L}^+ \cup \mc{L}^-)$, we simply argue by pseudo-convexity. Near $\Phi(\mc{L}^+ \cup \mc{L}^-)$, we use the construction of our ``barrier" in \Cref{lemma:hypersurface-separating}, along with a positivity of intersection argument.

\pf
We prove that the image of $\mc{U}^+$ under $\Psi^+$ is closed since the other case is analogous. By \Cref{lemma:unique-standard} the closed subset $\R \setminus (-M,M)$ is contained in the image. It is thus sufficient to show that any convergent sequence $s_j \to b$ of points $s_j \in (-M,M) \cap \Psi^+(\mc{U}^+)$ has a limit $b$ that again is contained inside the image. (Recall that $\Psi^{+}$ is an embedding, which makes it possible to identify points in the image with points in the domain.) 

The proof consists of two steps. First we show that any $u \in \mc{U}^{+}$ with $\Psi^{+}(u) \in [-M,M]$ can be confined a priori to some fixed bounded subset of $T^*\Sigma$ (under the inclusion $T^*\S - \Phi(\mc{L}^+ \cup \mc{L}^-) \hookrightarrow T^*\S$). Then we apply the SFT-compactness theorem to extract a limit configuration, and ultimately deduce closedness.

Choose a sequence $s_j \to b$ with $s_j \in (-M,M) \cap \Psi^+(\mc{U}^+)$ and consider the sequence of planes $u^{+}_{s_j}$ inside $T^*\Sigma - \Phi(\mc{L}^+ \cup \mc{L}^-).$ We start by showing that these planes are contained in a uniformly bounded domain (i.e. independent of $s_j$) inside $T^*\Sigma$. To that end, it is an important feature of the moduli space that it consists of planes are proper inside $T^*\Sigma - ( \mc{L}^+ \cup \mc{L}^-)$, and asymptotic to the family of simple closed geodesics on $\Phi(\mc{L}^+ \cup \mc{L}^-)$. We will treat separately boundedness in the fiber and base directions. 

\emph{Boundedness in the fiber direction:} Recall from \Cref{lemma:convex-hyp} that $$\mc{S}_{r}:= \{ (q, p) \in T^* \S \mid \langle p, p \rangle_{g_{\eta}} = r \} \sub T^*\S$$ is pseudoconvex with respect to $j_S$, for any $r>0$. Let $H: T^* \S \to \R$ be the function $H(q,p):= \langle p, p \rangle_{g_{\eta}}$. Choose $N>0$ so large that $\mc{K} \sub \{(q, p) \in T^*\S \mid H(q,p)< N \}$ and $N> 100 \op{max}(M, C)$ (recall that $\mc{K}$ is a compact subset of $T^*\S$ introduced in \Cref{subsection:geometric-setup}, while $C>0$ is a constant introduced in \Cref{subsection:geometric-setup} which is equal --up to a small error-- to the radius of the Lagrangian cylinders we have removed). Suppose that $\op{max} H \circ u^+_{s_j}= R \geq N^2$ and is achieved at some point $P=u^+_{s_j}(p)$, for $p$ a point in the domain of $u^+_{s_j}$.  

Let us first assume that $P$ is contained in the region $\{(x_1, y_1, x_2, y_2) \mid |y_2| <1/4 \}  \sub T^*\S_+$. We write $P= (P_{x_1}, P_{y_1}, P_{x_2}, P_{y_2})$.  According to \Cref{lemma:hypersurface-separating} and positivity of intersection, the curve $u^+_{s_j}$ cannot cross $\mc{H}$ and must therefore be contained entirely in $\mc{V}$. (Indeed, $\mc{H}$ is foliated by properly embedded, boundaryless $J$-holomorphic surfaces. If $u^+_{s_j}$ intersected $\mc{H}$, then by positivity of intersection, it must pair positively with one of the curves in the foliation. But this means that it pairs positively with all of the curves in the foliation, which is clearly false since these curves eventually become disjoint from $u^+_{s_j}$.) It follows that $|P_{x_2} | < 2M+1$, and hence $|P_{x_1}| > N^2-(3M)^2 > N^2/2$.\footnote{It follows from the construction of $g_{\eta}$ in \eqref{align:glued-metric} that the metric induces by $g_{\eta}$ on $T^*\S$ in the region $\{(x_1, y_1, x_2, y_2) \mid |y_2| <1/4 \}  \sub T^*\S_+$ coincides with the euclidean metric on $\R^4$ under the identification \eqref{equation:id1}. Thus there is no ambiguity in the meaning of $|P_{x_1}|$ and $|P_{x_2}|$.} Hence $J=j_S$ near $P$ by the definition of $J$ (see \cref{subsection:geometric-setup}, in particular note that $P \notin \mc{K}$ and that $P$ is far away from $\Phi(\mc{L}^\pm)$ (combine \eqref{equation:M-in-base} and the fact that $N> 100C$).  Hence $\mc{S}_R$ is $J$-convex near $P$, which gives a contradiction. 

If $P \in T^* \S $ contained in the complement of the region $\{(x_1, y_1, x_2, y_2) \mid |y_2| <1/4 \}$, then $J=j_S$ near $P$ by definition of $J$ (again $P \notin \mc{K}$, and $P$ is far away from $\Phi(\mc{L}^\pm)$). In particular, $\mc{S}_{R}$ is $J$-convex near $P$, which again gives a contradiction. We conclude that $u^+_{s_j}$ cannot cross the hypersurface $\mc{S}_{N^2}$, i.e. $u^+_{s_j}$ is bounded in the fiber directions. 

\emph{Boundedness in the base direction:} Let us first consider the function $\pi_{y_1} \circ u^+_{s_j}$, which is well defined on $(u^+_s)^{-1}(T^* \S_+)$. Observe that the hypersurfaces $\{y_1= N' \} \sub T^*\S$ are foliated by properly embedded, boundaryless, $J$ holomorphic planes if $|N'| \geq N$ (indeed, $J=j_S$ in that region). Hence, by positivity of intersection, $u^+_{s_j}$ cannot touch these hypersurfaces. It follows that $|\pi_{y_1} \circ u^+_{s_j}| < N'$.  A similar argument shows that $\pi_{y_2} \circ u^+_{s_j} < N$.

It remains to argue that $u^+_s$ cannot approach $B$ (see \eqref{equation:boundary}).  This is a consequence of the fact that $g_{\eta}= (\phi_k^{\eta})_*(g_e)$ in the region $\phi_k^{\eta}( \R / \Z \tms (0, \eta/4))$ for $1 \leq k \leq \beta$ (resp.\ in the region $\phi_0^{\eta}(\R \tms (0, \eta/4))$. Indeed, for $0<c< \eta/4$, we find that the properly embedded hypersurfaces $\pi^{-1}(\op{Im} \phi_k^{\eta}(-, c))$ are foliated by properly embedded $j_S$-holomorphic curves (in the case $k=0$ these are planes, while in the case $k>0$ these are cylinders). Hence the $u^+_{s_j}$ cannot touch them.

We conclude that the $u^{+}_{s_j}$ are contained in a uniformly bounded domain. Hence we can apply the SFT compactness theorem (see \Cref{subsection:punctured-curves}). 

For area reasons, the limit building consists of a single plane $v \in \mc{M}^+$. Moreover, since $u^{+}_{s_j} * u^{+}_{s_j}=0$, it follows that $v*v=0$.  Noting that $v$ must have vanishing normal Chern number (since the $u^{+}_{s_j}$ do), it follows from the adjunction formula for punctured holomorphic curves (see \Cref{fact:adjunction}) that $v$ is embedded. It then follows from automatic transversality (see \Cref{fact:automatic-transversality}) that $v$ is transversally cut out. Hence $v \in \mc{U}^{+}$ and $\Psi^+(v)=b$. \epf

By putting to together the previous lemmas, we find that $\Psi^{\pm}$ is a smooth embedding of $\mc{U}^{\pm}$ whose image is closed. It follows that $\Psi^{\pm}$ is surjective, and hence a diffeomorphism. We state this as a corollary. 
\cor
\label{cor:diffeo}
The map $\Psi^{\pm}$ is a diffeomorphism.
\qed
\ecor

Let us now recall some notation from \Cref{subsection:geometric-setup}. Recall from \eqref{equation:w'} that $W':= \bigcup_{l \in \R} W'_l$ is a hypersurface ``with holes", where $W'_l:= \Phi(W_l)$ and $W_l$ was defined in \eqref{equation:wl}. Noting that $W_l \sub X$ and $J= \Phi_*J'$ in $\Phi(X)$, one can check that $W'_l$ is $J$-holomorphic. Positivity of intersection and the invariance in families of the intersection number then implies the following properties which we collect as a lemma.

\lem \label{lemma:positivity}
The following properties hold:
\begin{itemize}
	\item[(i)] $W'_{l_1} \cap u^{\pm}_{l_2} =\emptyset $ for any $l_1, l_2 \in \R$.
	\item[(ii)] $u^+_{l_1} \cap u^-_{l_2}= \emptyset$ for any $l_1, l_2 \in \R$.
	\item[(iii)] $u^+_{l_1} \cap u^+_{l_2}= \begin{cases} u^+_{l_1}  &\text{if } l_1 = l_2, \\
	\emptyset &\text{otherwise}. \end{cases}$
	\item[(iv)] The analog of (iii) holds with ``$-$" in place of ``$+$". 
\end{itemize}
\qed
\elem 

We now fill-in the holes of $W'$. More precisely, the goal is to extend $W'$ to a hypersurface $Q \sub T^*\Sigma$ which is diffeomorphic to $\R^3$. Moreover, the symplectic leaves $W'_l \subset W'$ (which are each diffeomorphic to $\R^2$ with two disks removed) should extend to symplectic leaves $\S_l \sub Q$, where $\S_l$ is diffeomorphic to $\R^2$. 

To construct the $\S_l$, the idea is to adjoin to $W'_l$ one symplectic plane from each of the two moduli spaces $\mc{U}^{\pm}$. This does not quite work on the nose, because the result of performing such a gluing might not be smooth (the best we can say is that it is $C^0$). To get around this problem, we will argue as in \cite[Sec.\ 5.3]{dgi} and \cite[Lem.\ 5.5]{wendlduke}: there are asymptotic estimates which tell us that finite energy pseudoholomorphic planes converge to a trivial cylinder. Using these estimates, one can smoothly deform the family of planes so that they become set-wise \emph{equal} (i.e.\ not just asymptotic) to trivial cylinders, far enough inside the negative end. It is then immediate to show that the deformed family can be smoothly glued to $W'$ along the Lagrangian cylinders $\Phi(\mc{L}^+ \cup \mc{L}^-)$. 

\prop[Smoothing procedure]\label{proposition:smoothing-procedure}
\label{prop:smoothing}
Identify $\mc{U}^\pm \tms \C$ with $\R_l \tms \C$ via $\Psi^{\pm}: \mc{U}^\pm \xrightarrow{\sim} \R_l$. The corresponding family of parametrized pseudoholomorphic planes $u_l^{\pm} \colon \C \to T^*\S$ can be smoothly deformed in the concave end so that the following is satisfied: 
	\begin{itemize}
	\item[(i)] The sets $\S_l:=  \op{im}(u^{+}_l)  \cup W'_l \cup \op{im}(u^{-}_l)$ are codimension $2$ embedded symplectic surfaces. Moreover, they fit together to form a smooth foliation of their union $$Q:= \cup_l \S_l,$$ which is a properly embedded submanifold diffeomorphic to $\R^3$. 
	%\item[(ii)] For $l_1, l_2 \in \R$, we have $\S_{l_1} \cap \S_{l_2} = \begin{cases} \S_{l_1} &\text{if } l_1=l_2 \\
	%\emptyset &\text{otherwise}.  \end{cases} $
	\item[(ii)] For $|l|>M$, we have that $\S_l= \{x_2=l, y_2=0\} \sub T^*\S_+ \sub T^*\S$. 
	\item[(iii)] For any $l \in \R$, we have that $\S_l \cap \{\max\{|x_1|,|y_1|\}>M+1 \} = \{\max\{|x_1|,|y_1|\}>M+1, x_2=l, y_2=0 \}$.
	\end{itemize}
\eprop
Note that the identification \eqref{equation:id1} is implicit in (ii) and (iii) of the above proposition. Recall also that $W'$ contains $L$, so obviously $Q$ also contains $L$.
\pf
The entire argument takes place in the negative cylindrical ends of $T^*\S-\Phi(\mc{L}^+ \cup \mc{L}^-)$. Our first task is to review some notation from previous sections in order to extract nice coordinates for these negative ends. 

To begin with, we refer the reader back to \Cref{subsection:a-c-structures}. There, we considered $T^*(S^1_{\theta} \tms \R_r)$ with coordinates $(\t, t, r, s)$. We also considered (see \eqref{equation:sphere-bundle}) the sphere bundle $S_{\tilde{\e}}:= \{ (\t, t, r,s) \mid \|(t,s) \| = \tilde{\e} \} \sub T^*(S^1 \tms \R)$, which has coordinates $(\theta, r, \phi) \in S^1 \tms \R \tms S^1$. Finally, we described (see \eqref{equation:cyl-end}) a negative cylindrical end
\eq\label{equation:negative-ends} ((-\infty, 0]_\tau \tms S^*_{\tilde{\e}}, d(e^\tau \a), J_0). \eeq

We claim that the two negative ends of $(T^*\S-\Phi(\mc{L}^+ \cup \mc{L}^-), d\l_{\op{can}}, J)$ are isomorphic (both symplectically and as almost-complex manifolds) to \eqref{equation:negative-ends}. This identification is given by the composition $\Phi \circ F^{\pm} \circ \iota$, where the definition of $\iota$ is in \Cref{subsection:a-c-structures} and the definitions of $F^\pm, \Phi$ are in \Cref{subsection:geometric-setup}. We have already checked that these maps are indeed symplectic embeddings (see again \Cref{subsection:a-c-structures} and \Cref{subsection:geometric-setup}). To check that we also have an isomorphism of almost-complex manifolds, one needs to trace through the construction of $J$ (combine \eqref{equation:jcyl-j0} and \eqref{equation:cylindrical-ac-ref}). 

We also note that the zero section $0_{S^1 \tms \R} \sub T^*(S^1_{\theta} \tms \R_r)$ is identified with $\Phi(\mc{L}^\pm)$ via $\Phi \circ F^{\pm}$ (indeed, we already remarked in \Cref{subsection:geometric-setup} that $F^\pm$ takes the zero section to $\mc{L}^\pm$). 

The closed Reeb orbits on each component of $S^*_{\tilde{\e}}$ are in bijective correspondence with \emph{oriented} closed geodesics for the flat metric on $S^1_\t \tms \R_r$.
In fact, we can parametrize them as follows: $\eta(l): S^1 \ni \theta \mapsto (\theta, l, 0) \in  S^*_{\tilde{\e}}$ and $\nu(l): S^1 \ni \theta \mapsto (-\theta, l, \pi) \in  S^*_{\tilde{\e}}$. We also have trivial cylinders $c^{\eta}_l: (-\infty, 0] \tms S^1 \ni (\tau, \theta) \mapsto (\tau, \theta, l, 0) \in (-\infty, 0] \tms S^*_{\tilde{\e}}$ and $c^{\nu}_l: (-\infty, 0] \tms S^1 \ni (\tau, \theta) \mapsto (\tau, -\theta, l, \pi) \in (-\infty, 0] \tms S^*_{\tilde{\e}}$. 

We are going to use \eqref{equation:negative-ends} as our coordinate system to describe the negative ends of $T^*\S-\Phi(\mc{L}^+ \cup \mc{L}^-)$. In fact, we will focus entirely on the negative end corresponding to $\Phi(\mc{L}^+)$, which is identified with \eqref{equation:negative-ends} from now on (the other negative end can be treated similarly). %We therefore drop the subscripts $(-)^\pm$ in our notation when there is no ambiguity. 

We now make three important observations:
\begin{enumerate} 
\item The union of the trivial cylinders $c^{\eta}_l \cup c^{\nu}_l$ compactifies to form a family of smoothly embedded symplectic annuli. Moreover, these annuli intersect $0_{S^1 \tms \R} \equiv \Phi(\mc{L}^+)$ cleanly in the underlying closed geodesic $\{(\theta, l) \mid \theta \in S^1\}$. 
\item For all $l \in \R$, the trivial cylinders $c^{\eta}_l$ embed set-wise into $W'_l$. To check this, see \eqref{equation:wl}, \eqref{equation:def:w'l} and \eqref{align:def-of-F}. 
\item\label{item:linear-large-m} For $| l | >M$, the planes $u^+_l$ coincide with the standard ``linear" planes $u_{0, \s}^{+}$ (this follows from \Cref{lemma:unique-standard}; see also \eqref{equation:standar-linear-plane} for the definition of the ``linear" planes). As a result, the trivial cylinder $c^{\nu}_l$ embeds set-wise into $\op{im}(u^+_l)$ for $|l|>M$.
\end{enumerate}

Analogous statements hold in the negative end associated to $\Phi(\mc{L}^-)$, with $u_l^-$ in place of $u_l^+$. In particular, it follows from these observations that
$$ \op{im}(u^+_l)  \cup W'_l \cup \op{im}(u^-_l)$$
is a smoothly embedded symplectic submanifold of $T^*\S$ for $ | l | >M$. More generally, this holds whenever the planes $u^{\pm}_l$ coincide setwise with trivial cylinders outside a compact subset of their domain. %(See the proof of \cite[Prop.\ 5.16]{dgi} for a similar statement.) 

For $l \in [-M, M]$, we cannot in general guarantee that the $\op{im}(u^{\pm}_l)$ coincide setwise with a trivial cylinder outside a compact subset of their domain. However, following the arguments of \cite[Section 5.3]{dgi} or \cite[Lem.\ 5.5]{wendlduke}, one can deform these planes so that this property holds. We again explain this only for $u^+_l$. 

The key point is that these planes converge \emph{asymptotically} to trivial cylinders. Moreover, since the family of planes is standard ``linear" outside a compact set, the asymptotic convergence is uniform. More precisely, use the exponential map to identify a neighborhood of $\nu(l)$ with a neighborhood of the zero section in $\nu(l)^*\xi \to \nu(l)$. Then, up to smooth reparametrization, we can write $u_l^+(\tau, \theta)= (\tau, \s_l(\tau,\theta))$, where $\s_l(\tau, -)$ is a section of $\nu(l)^*\xi \to \nu(l)$. 

Now it can be shown (see e.g.\ \cite[Lem.\ 5.14]{dgi} or \cite[(A.3)]{wendlduke}\footnote{The constants $T, 2\pi, S_\l, t_\l$ in \cite[Lem.\ 5.14]{dgi} depend on how one chooses to parametrize trivial cylinders, and they can be absorbed by an appropriate reparametrization of the domain. In \cite[(A.3)]{wendlduke}, the asymptotic formula is stated in the Morse--Bott setting, but only for positive cylindrical ends. However, the version for negative ends is identical \emph{mutatis mutandis}.} ) that $(l,\theta) \mapsto (\s_l(\tau, \theta)- \nu(l)(\theta) )$ converges to zero in the $C^1$ norm as $\tau \to -\infty$, \emph{uniformly in $l$ and $\theta$}. (Here, we are of course viewing $\theta \mapsto \nu(l)(\theta)$ as the zero section of $\nu(l)^*\xi \to \nu(l)$.) The fact that the convergence is uniform in $l$ is the crucial point. The only difference between the argument in the present case and the argument in the aforementioned cases \cite{dgi} and \cite{wendlduke} is that our moduli space is not compact. However, since $\s_l(\tau,\theta) \equiv \nu(l)(\theta)$ for $| l | >M$ by \eqref{item:linear-large-m}, we can treat the moduli space similarly as in the compact case.  

It follows that it is possible to deform the family of planes $u^+_l$ with $|l| \le M$ to symplectic planes that coincide setwise with $c^{\nu}_l$ far enough inside the negative end (i.e.\ for $\tau \ll 0$). More precisely, let $\vartheta\colon (-\infty, 0] \to [0,1]$ be a non-decreasing smooth function which is $0$ on $(-\infty, -1/2]$ and $1$ on $[-1/4, 0]$. For $m\gg 0$ to be fixed later, let $\vartheta_m(-):= \vartheta(-/m)$. Now consider the interpolated family $u_l^{+,m}: (\tau, \theta) \mapsto (\tau, \vartheta_m(\tau) \s_l(\tau, \theta) + (1-\vartheta_m(\tau)) \nu(l)(\theta))$ (where it is understood that $u_l^{+,m}=u_l^+$ outside the negative end). 
% see.\ \cite[Lem.\ 5.16]{dgi} or \cite[Lem.\ 5.5]{wendlduke} for a similar argument.

We define $u^{-,m}_l$ similarly. We now set 
\eq \S_l:= \op{im}(u^{+,m}_l) \cup W'_l \cup \op{im}(u^{-,m}_l).\eeq

We claim that $\S_l$ satisfies the desired properties, provided that $m$ is large enough. We again only discuss the negative end associated to $\Phi(\mc{L}^+)$. To begin with: $\S_l$ is clearly smooth since $\op{im}(u^{+,m}_l)$ agrees set-wise with the trivial cylinder $c^{\nu}_l$ as one moves far enough into the negative end, so the components tautologically glue. The $C^1$ uniform convergence of $u_l$ to $c^{\nu}_l$ implies that the family is embedded for $m$ large enough (because the linear interpolation between $\theta \mapsto \nu(l)(\theta)$ and $\theta \mapsto \s_l(\tau,\theta)$ remains embedded). Finally, to see that the planes are symplectic, note that there is no difference in being symplectic with respect to $e^\tau(d \tau \wedge \a + d\a)$ or with respect to $(d \tau \wedge \a + d\a)$. The linear interpolation between $\theta \mapsto \nu(l)(\theta)$ and $\theta \mapsto \s_l(\tau,\theta)$ is uniformly $C^1$ close to a Reeb orbit for each $\tau$, so the desired claim follows easily for $m$ large enough. This establishes (i). The remaining properties (ii) and (iii) follow from the definition of $W_l'$ (see \eqref{equation:def:w'l}). 
\epf

For the remainder of this paper, the submanifolds $Q$ and $\S_l$ whose existence is guaranteed by \Cref{proposition:smoothing-procedure} are fixed.

\subsection{Completion of the proof} \label{subsection:completion-proof}
We start by recalling some general facts from symplectic geometry. Consider a hypersurface $\fk{H}^{2n-1}$ in a symplectic manifold $(M^{2n}, \o)$. There is a one-dimensional distribution $$\op{ker}\o|_{T\fk{H}} \sub T\fk{H}$$ called the \emph{characteristic distribution}. 

Suppose now that we have a symplectic fibration $\Pi: \fk{H} \to \R$ (i.e. a smooth fiber bundle whose fibers are symplectic with respect to the restriction of $\o$). Then the fibers are transverse to the characteristic distribution. Assuming that the fibration is suitably well-behaved at infinity (so that we can integrate horizontal vector fields), we can define a parallel transport map $\Pi^{-1}(l_0) \to \Pi^{-1}(l_1)$ for $l_1, l_2 \in \R$. It is a well-known fact, which can easily be checked using the Cartan formula, that parallel transport induces symplectomorphisms between the different fibers. 

Suppose now that $\L \sub \fk{H} \sub (M, \o)$ is a Lagrangian submanifold. Then
$$\Lambda \cap \Pi^{-1}(l_0) \subset (\Pi^{-1}(l_0),\o|_{T\Pi^{-1}(l_0)})$$
is itself a Lagrangian submanifold of one dimension lower than the dimension of $\Lambda$ (indeed, note that $\Pi^{-1}(l_0)$ and $\L$ must intersect transversally as submanifolds of $\fk{H}$). The intersections of the Lagrangian $\Lambda$ with the other fibers of $\Pi$ are then given as the image of $\Lambda \cap \Pi^{-1}(l_0)$ under parallel transport. Observe that the converse also holds, i.e.~ any Lagrangian submanifold of a given fiber $(\Pi^{-1}(l_0),\o_{can}|_{T\Pi^{-1}(l_0)})$ gives rise to a Lagrangian submanifold of $(M, \o)$ contained in $\fk{H}$ by extending it to a submanifold which is tangent to the characteristic distribution.

If $(M^4, \o)$ is a symplectic $4$-manifold, then any real one-dimensional curve in $\Pi^{-1}(l_0)$ is automatically Lagrangian. This makes Lagrangians in $\fk{H} \sub M^4$ particularly easy to understand and construct: indeed, by the argument above, they are determined by a one-dimensional curve inside any one of the fibers.

Let us now specialize the above discussion to the case which is relevant from proving \Cref{theorem:main}. By modifying the hypersurface $Q \sub T^*\S$ obtained in \Cref{prop:smoothing}, we will construct a hypersurface $\tilde{Q} \sub T^*\S$ satisfying the following three properties. 

\begin{enumerate}
\item[(P1)] \label{enumerate:standard-infinity}  There is a smooth symplectic fibration $\tilde{Q} \to \R$ whose symplectic fibers $\tilde{\S}_l$ are proper embeddings $\R^2 \hookrightarrow T^*\S$, and such that the fiber over $l \in \R$ coincides with the symplectic plane
\eq \label{equation:symplectic-planes} \{x_2=l,y_2=0\} \subset T^*\S \eeq
outside of a compact subset. In addition, there exists $N \geq M>0$ (where $M$ was defined in \Cref{section:standard-model}) so that the fibers for $|l| \geq N$ coincide with the planes \eqref{equation:symplectic-planes}. 
\item[(P2)]  \label{enumerate:contains-L} The Lagrangian $L$ is contained in $\tilde{Q}$ (and consequently the characteristic distribution of $\tilde{Q}$ is tangent to $L$).
\item[(P3)]  \label{enumerate:foliation-standard} The Lagrangian submanifolds $\chi_s \subset \tilde{Q}$ which are uniquely determined by the requirement that
$$\chi_s \cap \tilde{\S}_{-N}=\{y_1=s,x_2=-N,y_2=0\}=: \ell_s$$
satisfy $$\chi_s \cap \tilde{\S}_N=\{y_1=s,x_2=N,y_2=0\}$$ for any $s \ge 0$.  
\end{enumerate}

It follows from the above discussion that the $\chi_s$ are diffeomorphic to $\R^2$ and that $\chi_0=L$ (indeed, $\chi_0 \cap \tilde{\S}_{-N}= \ell_0 =L \cap \tilde{\S}_{-N}$, since $N\geq M$). 

An immediate consequence of (P1) is that the characteristic distribution of $\tilde{Q}$ is equal to $\R \cdot \partial_{x_2}$ outside of a compact subset. The parallel transport from the fiber over $-N$ to the fiber of $N$ will be called the \emph{symplectic monodromy map} $$\mu: \tilde{\S}_{-N} \to \tilde{\S}_N.$$ Note that it follows from (P1) that $\mu$ has compact support. Since $L \cap \tilde{\S}_{-N} = \ell_0=\{y_1=0\} \cap \tilde{\S}_{-N}$, it follows from (P1) and (P2) that $\mu(\tilde{\S}_{-N} \cap \{y_1=0\})=\tilde{\S}_{N} \cap \{y_1=0\} =\tilde{\S}_N \cap L$. Using now (P3), the sought Lagrangian isotopy that unknots $L$ is easy to construct: all that one has to do is to translate the Lagrangian planes $\chi_s$ appropriately in order to make the family compactly supported. The detailed argument will be given in \Cref{corollary:final}. We remark that the general technique for constructing isotopies of Lagrangians using the characteristic distribution of a hypersurface already appears in \cite{eli-pol2}.

We now explain how to construct $\tilde{Q}$ from $Q$.  Observe that \Cref{prop:smoothing} furnishes an obvious symplectic fibration $Q \to \R$ with fibers $\S_l$. The fibers are standard at infinity according to \Cref{prop:smoothing}, so we have a well-defined parallel transport map $\S_l \to \S_{l'}$ for $l, l' \in \R$. Let us now set $M':= M+2$. 

Let $\mu^+$ denote the restriction the parallel transport map $\S_{-M'} \cap \{y_1 \geq 0\} \to \S_{M'} \cap \{y_1 \geq 0\}$ (observe that this is well defined since $L \sub Q$ and $L$ is preserved by the characteristic flow). According to \cite[Lem.\ 6.8]{dgi}, $\mu^+$ is generated by a family of compactly supported Hamiltonians $\{H_t\}_{t \in [0,1]}$ on $\{(x_1, y_1) \mid y_1 \geq 0\}$ which vanish on $\{y_1=0\}$ and such that $H_t \equiv 0$ for $t$ near $\{0,1\}$.  We will use this Hamiltonian to deform $Q$ via the ``symplectic suspension'' construction, that we now outline.

Given $R>0$, let $H_t^{R}:= \frac{1}{R} H_{Rt}$. Let $\phi^{H^{R}}_{t}$ be the associated Hamiltonian flow where $t \in [0,R]$. Set $x_2':= x_2- M'$. 

We now fix $R$ large enough so that $| H_t^R|< 1/2$.

\lem \label{lemma:symplectic-embedding}
We have a symplectic embedding $$\Theta: T^* \S_+ \cap \{ M' \leq x_2 \leq M'+R, y_2> -1/2 \} \to T^*\S_+ \sub T^* \S$$ given by setting $\Theta(x_1, y_1, x_2, y_2)= (\phi^{{H}^{R}}_{x_2'}(x_1, y_1), x_2, y_2+ H^{R}_{x_2'})$. 
\elem
\pf
Note first of all that this map is well-defined due to our assumption that $|H_t^R|<1/2$ and the definition of $\S_+ \sub \S$ in \Cref{section:standard-model}. 

We let $\phi_{x'_2}^x= \pi_{x_1} \circ \phi^{H^R}_{x'_2}$ and let $\phi_{x'_2}^y= \pi_{y_1} \circ \phi^{H^R}_{x'_2}$. We also denote differentiation in $x'_2$ by a dot. 

We now compute
\begin{align*} 
\Theta^*\o &=  \d_{x_1} \phi^x dx_1 \wedge dy_1+ \d_{y_1}\phi^y dx_1 \wedge dy_1 + \dot{\phi}^x dx_2 \wedge dy_1 + \dot{\phi}^y dx_1 \wedge dx_2 + dx_2 \wedge dy_2+ \\&\qquad dx_2 \wedge (\d_{x_1}H_{x'_2}^R dx_1 + \d_{y_1} H_{x'_2}^R dy_1) \\
&=(\d_{x_1} \phi^x dx_1 \wedge dy_1+ \d_{y_1}\phi^y dx_1 \wedge dy_1 + dx_2 \wedge dy_2) + \dot{\phi}^x dx_2 \wedge dy_1 + \dot{\phi}^y dx_1 \wedge dx_2 +\\&\qquad dx_2 \wedge (\dot{\phi}^y dx_1 - \dot{\phi}^x dy_1) \\
&=\o.
\end{align*}
\epf

Now let
$$
\tilde{Q} = \begin{cases}  
{Q} &\text{for } x_2 \in \R - [M', M'+R], \\
\op{Im} \Theta({Q} \cap \{ N \leq x_2 \leq M'+R \}) &\text{for } x_2 \in [M', M'+R].
\end{cases}
$$

Since $H_t \equiv 0$ for $t$ near $\{0,1\}$, it follows that $\tilde{Q}$ is a smooth hypersurface which agrees with $Q$ outside a compact set. 

\begin{proof}[Proof that $\tilde{Q}$ satisfies (P1)--(P3)]
First of all, \Cref{prop:smoothing} implies that (P1) and (P2) are satisfied by $Q$. Since the Hamiltonian $H_t$ has compact support, it follows from the construction that these properties are still satisfied by $\tilde{Q}$ with $N= M'+R$.

Property (P3) is a consequence of the stronger claim that the monodromy map $\mu$ for $\tilde{Q}$ restricts to the identity on the set $\{y_1 \geq 0, x_2=-N, y_2=0\} \sub \tilde{\S}_{-N}$. This stronger claim is in turn a consequence of \Cref{lemma:symplectic-embedding}. Indeed, the fact that $\Theta$ is a symplectomorphism implies that the characteristic distribution of $Q$ is mapped to the characteristic distribution of $\tilde{Q}$ under $\Theta$; in other words, the monodromy map of $\tilde{Q}$ is obtained from that of $Q$ by simply post-composing with the restriction $\Theta|_{\Sigma_{M'+R}}$.
\end{proof}

Having proved that $\tilde{Q}$ satisfies (P1)--(P3), we obtain the following corollary, which immediately implies \Cref{theorem:main} (recall from \Cref{section:standard-model} that we may assume without loss of generality that $x=0$ in the statement of \Cref{theorem:main}). 

\cor \label{corollary:final}
There exists a compactly supported Hamiltonian isotopy of $(T^*\S, d\l_{can})$ taking $L$ to $F_0$.
\ecor

\pf
For $s \geq 0$, recall that $\ell_{s} \sub \tilde{Q}$ is the line $\ell_{s}:= \{y_1=s, x_2= -3N, y_2= 0\}$. We also have that $\chi_s \sub \tilde{Q}$ is a Lagrangian plane which is standard at infinity and is obtained as the image $\ell_{s}$ under the characteristic flow. Finally, we have that $\chi_0=L$ and that $\chi_{s_0}= \{y_1=s_0, y_2= 0\} \sub \tilde{Q}$ for $s_0$ large enough. 
	
Let $\s^+: [0, s_0] \tms \S_+ \to \S_+$ be a compactly supported isotopy with the property that $\s_s^+(s, 0)=(0, 0)$ and that $\s_s(y_1, \a)= (y_1, \a)$ if $ \a \in [-1, -1/2]$. Observe that $\s^+$ extends to an isotopy $\s : [0, s_0] \tms \S \to \S$ which is the identity on $\S- \S_+$. 
	
We now consider the induced isotopy $\s^*: [0, s_0] \tms T^* \S \to T^*\S$. Observe that $s \mapsto \s_s^*(\chi_s)$ defines a smooth family of Lagrangian submanifolds which are diffeomorphic to $\R^2$ and fixed setwise outside a compact set. Choose a compactly supported smooth isotopy of $T^*\S$ which gives rise to the family $\s_s^*(\chi_s)$ of embedded submanifolds. The infinitessimal generator $\zeta_t\in\Gamma(TT^*\S)$ of the ambient isotopy satisfies the property that
$$\alpha_s \coloneqq -(\iota_{\zeta_s}\omega_{can})|_{T\s_s^*(\chi_s)}$$
is a smooth family of compactly supported one-forms on $\s_s^*(\chi_s)$. By Cartan's formula, we have $d\alpha_s= -(\mc{L}_{\zeta_s} \o)|_{T\s_s^*(\chi_s)}=0$, where the last equality follows since $\s_s^*(\chi_s)$ is Lagrangian. Thus the $\alpha_s$ are closed. Since $H^1_c(\R^2; \R)=0$ it follows that there is a smooth family of compactly supported smooth functions
$$H_t \colon \s_s^*(\chi_s) \to \R$$
for which $dH_t=\alpha_t$. An arbitrary compactly supported extension of $H_t$ to $T^*\S$ now gives rise to the desired Hamiltonian isotopy.
\epf

\end{document}